# Derivation and Analysis of the Primal-Dual Method of Multipliers Based on Monotone Operator Theory

Thomas Sherson, Richard Heusdens, and W. Bastiaan Kleijn

*Abstract*—In this paper we present a novel derivation for an existing node-based algorithm for distributed optimisation termed the primal-dual method of multipliers (PDMM). In contrast to its initial derivation, in this work monotone operator theory is used to connect PDMM with other first-order methods such as Douglas-Rachford splitting and the alternating direction method of multipliers thus providing insight to the operation of the scheme. In particular, we show how PDMM combines a lifted dual form in conjunction with Peaceman-Rachford splitting to remove the need for collaboration between nodes per iteration. We demonstrate sufficient conditions for strong primal convergence for a general class of functions while under the assumption of strong convexity and functional smoothness, we also introduce a primal geometric convergence bound. Finally we introduce a distributed method of parameter selection in the geometric convergent case, requiring only finite transmissions to implement regardless of network topology.

## I. INTRODUCTION

Over the last three decades, there has been a significant rise in the interest in and deployment of large scale sensor networks for a wide range of applications [1], [2], [3]. Motivated by the increase in computational power of low cost micro-processors, the target applications for such sensing networks have become more diverse and computationally demanding. Applications of such networks include environmental monitoring [4], [5], power grid management [6], [7], [8], as well as part of smart homes enabling the likes of home health care systems [9], [10]. Additionally, where centralised network topologies were once the port of call for handling data processing of such systems, increasingly on-node computational capabilities of such networks are exploited to parallelise or even fully distribute data processing and computation. In contrast to their centralised counterparts such *distributed networks* have a number of distinct advantages including robustness to node failure, scalability with network size and localised transmission requirements.

Distributed networks, characterised by their limited connectivity, implicitly restrict the data available at any one node. As this data is often generated within the network, this limited connectivity, in combination with the removal of a single point of computation, makes even straightforward signal processing, such as the computation of inner products, cumbersome. To address this point, within the literature, a number of algorithms have been developed to solve a variety of signal processing problems in a distributed context. In the case of linear averaging problems, methods such as randomised gossip [11], weighted gossip [12] and path averaging [13] have received much attention due to their simplicity of implementation. However, the limited scope of these algorithms restricts their use in more general distributed settings.

Another family of distributed algorithms are those based on message-passing or belief propagation schemes [14], [15]. Including the sum-product and min-sum algorithms, such methods can be used to compute marginal distributions and solve maximum a posteriori (MAP) problems in a distributed fashion. By embedding a desired operation within such MAP problems, these methods can be used for more general distributed computation than averaging alone. Unfortunately, such methods are only provably convergent for acyclic networks, which may not be present in practical applications.

A third approach commonly used in the literature, in applications such as sensor localisation [16], global averaging of data [12] and network utility optimisation [17], is to embed desired signal processing operations in a convex optimisation problem, which can then be solved in a distributed manner. Commonly used distributed convex solvers leverage sub-gradient [18], [19], and primal-dual algorithms [20], where each algorithm is chosen based on the specifications of a particular application. Solvers based on classical algorithms such as proximal gradient [21], Douglas-Rachford splitting [22] and the alternating direction method of multipliers (ADMM) [23], have also seen a resurgence in the literature due to their applicability to distributed optimisation. In addition, the broad class of problems that can be transformed into convex programs makes these methods highly attractive.

In more recent years, a number of novel approaches for distributed optimisation, both convex and non-convex, have also been proposed. In the convex case, the works of [24], [25], [26], echoing advances in three term operator splitting such as Vu-Condat splitting [27], [28], provide general frameworks for distributed convex optimisation. Including more classical approaches, such as ADMM, as special cases, these algorithms leverage primal-dual type schemes to exploit functional separability to create fully distributed implementations. In contrast, the work in [29], [30] focuses on the more general problem of potentially non-convex optimisation. In particular, by at

Thomas Sherson is with the Department of Microelectronics, Circuits and Systems group, Delft University of Technology, The Netherlands. Email: t.sherson@tudelft.nl

Richard Heusdens is with the Department of Microelectronics, Circuits and Systems group, Delft University of Technology, The Netherlands. Email: r.heusdens@tudelft.nl

W. Bastiaan Kleijn is with the Department of Microelectronics, Circuits and Systems group, Delft University of Technology, The Netherlands. and with the School of Engineering and Computer Science, Victoria University of Wellington, New Zealand. Email: w.b.kleijn@tudelft.nl

each iteration approximating both objective and constraints with specific strongly convex and smooth surrogates, the proposed methods have provable guarantees on convergence to local minima. Furthermore, in contrast to other methods, the proposed approach need not explicitly require functional separability, only the separability of the surrogates used. This allows for the optimisation of problems typically outside of the scope of distributed algorithms.

Recently, another new algorithm for distributed optimisation was proposed within the literature termed the primal dual method of multipliers (PDMM) [31]. While demonstrating promising performance, including guaranteed average convergence, which in some example applications was faster than other methods such as ADMM [31], [32], in addition to an asynchronous method of implementation, there are still a number of open questions surrounding the approach. The purpose of this paper is therefore two-fold. Firstly we introduce a new derivation for PDMM, constructed from the perspective of monotone operator theory. In contrast to its original derivation, this approach allows PDMM to be linked with classical first order methods within the literature such as Peaceman-Rachford and Douglas-Rachford splitting techniques as well as the alternating direction method of multipliers (ADMM). Secondly, we use this approach to demonstrate new and stronger convergence results for PDMM. The most important of these points is its geometric convergence under strong convexity and smoothness assumptions. This point was demonstrated for ADMM in [33] but, prior to this paper, no such results exist for PDMM. In this way, this work also strengthens the performance guarantees of PDMM, an important point for practical distributed optimisation.

The remainder of this paper is organised as follows. In Sec. II and III we introduce appropriate nomenclature and background definitions to support the remainder of the document. In Sec. IV we introduce a monotone operator based derivation of PDMM based on a specific dual lifting approach. Sec. V demonstrates a range of convergence results for the PDMM algorithm under different assumptions of the distributed optimisation problems including sufficient conditions for primal optimality. This is strengthened in Sec. VI were we demonstrate primal geometric convergence under the assumption of strong convexity and functional smoothness. We also highlight a distributed method for parameter selection in this case. Finally, Sec. VII includes simulation results to reinforce and verify the underlying claims of the document and the final conclusions are drawn in Sec. VIII

## II. NOMENCLATURE

In this work we will denote by $\mathbb{R}$ the set of real numbers, by $\mathbb{R}^N$ the set of real column vectors of length $N$ and by $\mathbb{R}^{M \times N}$ the set of $M$ by $N$ real matrices. Let $\mathcal{X}, \mathcal{Y} \subseteq \mathbb{R}^N$. A set valued operator $\mathbf{T} : \mathcal{X} \to 2^{\mathcal{Y}}$ is defined by its graph, $\mathrm{gra}\,(\mathbf{T}) = \{(\mathbf{x},\mathbf{y}) \in \mathcal{X} \times \mathcal{Y} \mid \mathbf{y} \in \mathbf{T}(\mathbf{x})\}$ where $2^{\mathcal{Y}}$ is the power set of the set $\mathcal{Y}$. Similarly, the notion of an inverse of an operator $\mathbf{T}^{-1}$ is also defined via its graph so that $\mathrm{gra}\,(\mathbf{T}^{-1}) = \{(\mathbf{y},\mathbf{x}) \in \mathcal{Y} \times \mathcal{X} \mid \mathbf{y} \in \mathbf{T}(\mathbf{x})\}$. We denote by $\mathbf{J}_{\mathbf{T},\rho} = (\mathbf{I} + \rho \mathbf{T})^{-1}$ the resolvent of an operator and by $\mathbf{R}_{\mathbf{T},\rho} = 2\mathbf{J}_{\mathbf{T},\rho} - \mathbf{I}$ the reflected resolvent. The composition of two operators $\mathbf{T} : \mathcal{X} \to 2^{\mathcal{Y}}$ and $\mathbf{S} : \mathcal{Y} \to 2^{\mathcal{Z}}$ is given by $\mathbf{S} \circ \mathbf{T} : \mathcal{X} \to 2^{\mathbf{Z}}$. The set of fixed-points of $\mathbf{T}$ is denoted by $\mathrm{fix}\,(\mathbf{T}) = \{\mathbf{x} \in \mathcal{X} \mid \mathbf{T}(\mathbf{x}) = \mathbf{x}\}$. We denote by $|\bullet|$ the cardinality of a set. If $|\mathbf{T}(\mathbf{x})| = 1 \; \forall \; \mathbf{x} \in \mathcal{X}$ then $\mathbf{T} = f$ is a single valued. In this case $f^*$ denotes the Fenchel conjugate of $f$ defined as $f^*(\mathbf{y}) \triangleq \sup_{\mathbf{x}}\,(\mathbf{y}^T\mathbf{x} - f(\mathbf{x}))$ and $\partial f$ denotes the subdifferential of $f$. Furthermore, the set of all closed, convex, and proper (CCP) functions $f : \mathcal{H}^N \to \mathbb{R} \cup \{+\infty\}$ is denoted by $\Gamma_0\,(\mathcal{H}^N)$. Finally, if $\mathbf{T} : \mathcal{H}^N \to \mathcal{H}^M$ is a linear operator, we denote by $\mathrm{ran}\,(\mathbf{T})$ and $\mathrm{ker}\,(\mathbf{T})$ the range and kernel of $\mathbf{T}$ respectively.

## III. BACKGROUND

We begin by introducing some basic definitions that will be used throughout this paper (see [34] for a detailed overview).

**Definition III.1.** *Monotonicity: An operator $\mathbf{T} : \mathcal{H} \to 2^{\mathcal{H}}$ is monotone iff*

$$\langle \mathbf{u} - \mathbf{v}, \mathbf{x} - \mathbf{y} \rangle \geq 0 \quad \forall (\mathbf{x},\mathbf{u}), (\mathbf{y},\mathbf{v}) \in \mathrm{gra}\,(\mathbf{T}),$$

*is strictly monotone if*

$$\langle \mathbf{u} - \mathbf{v}, \mathbf{x} - \mathbf{y} \rangle > 0 \quad \forall (\mathbf{x},\mathbf{u}), (\mathbf{y},\mathbf{v}) \in \mathrm{gra}\,(\mathbf{T}), \mathbf{x} \neq \mathbf{y},$$

*and is $\sigma$-strongly monotone if*

$$\langle \mathbf{u} - \mathbf{v}, \mathbf{x} - \mathbf{y} \rangle \geq \sigma \|\mathbf{x} - \mathbf{y}\|^2 \quad \forall (\mathbf{x},\mathbf{u}), (\mathbf{y},\mathbf{v}) \in \mathrm{gra}\,(\mathbf{T}).$$

*Furthermore, an operator $\mathbf{T}$ is maximal monotone if*

$$\forall (\hat{\mathbf{x}}, \hat{\mathbf{u}}) \notin \mathrm{gra}(\mathbf{T}) \; \exists \; (\mathbf{x}, \mathbf{u}) \in \mathrm{gra}(\mathbf{T}) \mid \langle \mathbf{u} - \hat{\mathbf{u}}, \mathbf{x} - \hat{\mathbf{x}} \rangle < 0.$$

**Definition III.2.** *Nonexpansiveness and Lipschitz: An operator $\mathbf{T} : \mathcal{H} \to 2^{\mathcal{H}}$ is nonexpansive if*

$$\|\mathbf{u} - \mathbf{v}\| \leq \|\mathbf{x} - \mathbf{y}\| \quad (\mathbf{x},\mathbf{u}), (\mathbf{y},\mathbf{v}) \in \mathrm{gra}\,(\mathbf{T}),$$

*and strictly nonexpansive if*

$$\|\mathbf{u} - \mathbf{v}\| < \|\mathbf{x} - \mathbf{y}\| \quad \forall (\mathbf{x},\mathbf{u}), (\mathbf{y},\mathbf{v}) \in \mathrm{gra}\,(\mathbf{T}), \mathbf{x} \neq \mathbf{y},$$

*Furthermore, $\mathbf{T}$ is L-Lipschitz if*

$$\|\mathbf{u} - \mathbf{v}\| \leq L\|\mathbf{x} - \mathbf{y}\| \quad \forall (\mathbf{x},\mathbf{u}), (\mathbf{y},\mathbf{v}) \in \mathrm{gra}\,(\mathbf{T}).$$

*and is contractive if $L < 1$.*

**Definition III.3.** *Convexity: A function $f$ is convex if*

$$f(\theta\mathbf{x} + (1-\theta)\mathbf{y}) \leq \theta f(\mathbf{x}) + (1-\theta)f(\mathbf{y}) \quad \forall \mathbf{x}, \mathbf{y} \in dom\,(f)$$

*is strictly convex if*

$$f(\theta\mathbf{x} + (1-\theta)\mathbf{y}) < \theta f(\mathbf{x}) + (1-\theta)f(\mathbf{y})$$
$$\forall \mathbf{x}, \mathbf{y} \in dom\,(f), \mathbf{x} \neq \mathbf{y},$$

*and is $\sigma$-strongly convex if $f - \frac{\sigma}{2}\|\bullet\|^2$ is convex.*

It follows that $\partial f$ is $\sigma$-strongly monotone if $f$ is $\sigma$-strongly convex. Furthermore, $\partial f$ is maximal if $f \in \Gamma_0\,(\mathcal{H})$.

**Definition III.4.** *Smoothness: A function $f$ is L-smooth if it is differentiable and $\frac{L}{2}\|\bullet\|^2 - f$ is convex.*

Equivalently it can be shown that an $L$-smooth function satisfies the inequality

$$f(\mathbf{x}) \leq f(\mathbf{y}) + \langle \nabla f(\mathbf{y}), \mathbf{x} - \mathbf{y} \rangle + \frac{L}{2}\|\mathbf{x}-\mathbf{y}\|^2 \quad \forall \mathbf{x}, \mathbf{y} \in \mathrm{dom}(f).$$

Thus the gradient of an $L$-smooth function is $L$-Lipschitz.

## IV. A Monotone Derivation of the Primal-Dual Method of Multipliers

In this section we reintroduce a node based algorithm for distributed optimisation termed the Primal-Dual method of multipliers (PDMM) [31]. Unlike earlier efforts within the literature [31], [32], here we demonstrate how PDMM can be derived from the perspective of monotone operator theory. Additionally, a connection previously unnoted in the literature between PDMM and distributed ADMM is highlighted.

### A. Node Based Distributed Optimisation

Consider a general undirected network consisting of $N$ nodes. The associated graphical model of such a network is given by $\mathbf{G}(V, E)$ where $V = \{1, ..., N\}$ and $E$ denote the set of nodes and undirected edges respectively. Assume that each node is equipped with a function $f_i \in \Gamma_0\left(\mathbb{R}^{M_i}\right) \ \forall i \in V$ parameterised by a local variable $\mathbf{x}_i \in \mathbb{R}^{M_i}$. Under this model, consider solving the following optimisation problem in a distributed manner.

$$\min_{\mathbf{x}_i \ \forall i \in V} \sum_{i \in V} f_i(\mathbf{x}_i) \quad (1)$$
$$\text{s.t} \quad \mathbf{A}_{i|j}\mathbf{x}_i + \mathbf{A}_{j|i}\mathbf{x}_j = \mathbf{b}_{ij} \quad \forall (i,j) \in E.$$

Here $\mathbf{A}_{i|j} \in \mathbb{R}^{M_{ij} \times M_i}$ and $\mathbf{b}_{ij} \in \mathbb{R}^{M_{ij}}$. Note the distinction between the subscripts $i|j$ and $ij$. The prior is a *directional* identifier used to denote the directed edge from node $i$ to node $j$, while the later is an *undirected* identifier. Furthermore, let $M_V = \sum_{i \in V} M_i$ and $M_E = \sum_{(i,j) \in E} M_{ij}$.

We will refer to these problems as distributed convex optimisation problems where $\mathbf{A}_{i|j}$ and $\mathbf{b}_{ij}$ impose affine constraints between the local variables at neighbouring nodes.

The prototype problem in (1) includes, as a subset, the family of distributed consensus problems which minimise the sum of the local cost functions under network wide consensus constraints. Such problems have received much attention in the literature in recent years [11], [23], [35], [36]. The algorithm presented in this paper can therefore be used for this prupose.

The associated Lagrange dual [37] of (1) is given by

$$\min_{\boldsymbol{\nu}} \sum_{i \in V} \left( f_i^* \left( \sum_{j \in \mathcal{N}(i)} \mathbf{A}_{i|j}^T \boldsymbol{\nu}_{ij} \right) - \sum_{j \in \mathcal{N}(i)} \frac{\mathbf{b}_{ij}^T}{2} \boldsymbol{\nu}_{ij} \right), \quad (2)$$

where each $\boldsymbol{\nu}_{ij} \in \mathbb{R}^{M_{ij}}$ denotes the dual vector variable associated with the constraint function along edge $(i,j)$, $\mathcal{N}(i) = \{j \in V \mid (i,j) \in E\}$ is the set of neighbours of node $i$ and $f_i^*$ is the Fenchel conjugate of $f_i$. By inspection, each $\boldsymbol{\nu}_{ij}$ in (2) is utilised in two conjugate functions, $f_i^*$ and $f_j^*$ resulting in a coupling between neighbouring nodes.

To address the linking of the objective terms, the dimension of the dual problem can be lifted by introducing additional directed edge variables $\boldsymbol{\lambda}_{i|j}$ and $\boldsymbol{\lambda}_{j|i}$. These are then constrained so that at optimality $\boldsymbol{\lambda}_{i|j} = \boldsymbol{\lambda}_{j|i}$, leading to what we term as the *extended dual* of Eq (1)

$$\min_{\boldsymbol{\lambda}} \sum_{i \in V} \left( f_i^* \left( \sum_{j \in \mathcal{N}(i)} \mathbf{A}_{i|j}^T \boldsymbol{\lambda}_{i|j} \right) - \sum_{j \in \mathcal{N}(i)} \frac{\mathbf{b}_{ij}^T}{2} \boldsymbol{\lambda}_{i|j} \right)$$
$$\text{s.t.} \quad \boldsymbol{\lambda}_{i|j} = \boldsymbol{\lambda}_{j|i} \quad \forall i \in V, j \in \mathcal{N}(i). \quad (3)$$

While seemingly straightforward, it is this choice of lifting in the dual domain that distinguishes our contribution from existing approaches.

**Remark 1.** *In the specific case that all $\mathbf{A}_{i|j} = \mathbf{I}$, the extended dual form reflects a similar structure to that demonstrated in (Sec. 3.3) [38] where by the duplication of the edge variables was also considered. The method in [38] focused on the case of stochastic coordinate descent based approaches. In contrast, our motivation for the use of the extended dual form here, which considers the case of arbitrary $\mathbf{A}_{i|j}$ matrices, is to allow us to link PDMM with existing algorithms within the literature via monotone operator theory.*

We denote by $\boldsymbol{\lambda}$ the stacked vector of all $\boldsymbol{\lambda}_{i|j}$. The ordering of this stacking is given by $1|2 < 1|3 < \cdots < 1|N < 2|1 < 2|3 < \cdots < N|N-1$. In the case of a fully connected network the vector $\boldsymbol{\lambda} \in \mathbb{R}^{M_E}$ is given by

$$\boldsymbol{\lambda} = [\boldsymbol{\lambda}_{1|2}; \cdots ; \boldsymbol{\lambda}_{1|N}; \boldsymbol{\lambda}_{2|1}; \cdots, \boldsymbol{\lambda}_{N|N-1}].$$

For the primal constraints we introduce the matrix $\mathbf{C} \in \mathbb{R}^{M_E \times M_V}$ and vector $\mathbf{d} \in \mathbb{R}^{M_E}$ where in the case of a fully connected network

$$\mathbf{C} = \begin{bmatrix} \mathbf{C}_1 & \cdots & \mathbf{0} \\ \vdots & \ddots & \vdots \\ \mathbf{0} & \cdots & \mathbf{C}_N \end{bmatrix}$$
$$\mathbf{C}_i = [\mathbf{A}_{i|1}; \mathbf{A}_{i|2}; \cdots ; \mathbf{A}_{i|i-1}; \mathbf{A}_{i|i+1}; \cdots ; \mathbf{A}_{i|N}] \ \forall i \in V$$
$$\mathbf{d} = [\mathbf{d}_1; \mathbf{d}_2; \cdots ; \mathbf{d}_N]$$
$$\mathbf{d}_i = \frac{1}{2}[\mathbf{b}_{i1}; \mathbf{b}_{i2}; \cdots ; \mathbf{b}_{i(i-1)}; \mathbf{b}_{i(i+1)}; \mathbf{b}_{i2}; \cdots ; \mathbf{b}_{iN}] \ \forall i \in V. \quad (4)$$

For other network topologies the unnecessary rows of $\boldsymbol{\lambda}$, $\mathbf{C}$ and $\mathbf{d}$, corresponding to non-existent edges, can be removed.

We further introduce the symmetric permutation matrix $\mathbf{P}$ that maps between each pair of variables $\boldsymbol{\lambda}_{i|j}$ and $\boldsymbol{\lambda}_{j|i}$ so that the constraints in (3) can be rewritten as $\boldsymbol{\lambda} = \mathbf{P}\boldsymbol{\lambda}$. Finally, we define the function

$$f(\mathbf{x}) : \mathbb{R}^{M_V} \mapsto \mathbb{R} \triangleq \sum_{i \in V} f_i(\mathbf{x}_i)$$

as the sum of all local functions where $\mathbb{R}^{M_V} = \mathbb{R}^{M_1} \times \mathbb{R}^{M_2} \times \ldots \times \mathbb{R}^{M_N}$. It follows that (3) is equivalent to

$$\min_{\boldsymbol{\lambda}} f^*(\mathbf{C}^T \boldsymbol{\lambda}) - \mathbf{d}^T \boldsymbol{\lambda} + \iota_{\mathrm{ker}(\mathbf{I}-\mathbf{P})}(\boldsymbol{\lambda}), \quad (5)$$

where $\iota_{\mathrm{ker}(\mathbf{I}-\mathbf{P})}$ is the indicator function of the edge based constraints defined as

$$\iota_{\mathrm{ker}(\mathbf{I}-\mathbf{P})}(\mathbf{y}) = \begin{cases} 0 & (\mathbf{I}-\mathbf{P})\mathbf{y} = 0 \\ +\infty & \text{otherwise}. \end{cases}$$

*B. From an Unconstrained Optimisation Problem to a Non-expansive PDMM Operator*

From the unconstrained optimisation problem in (5) we now introduce a method of solving (1) in a distributed manner. We do so via the relationship between (5) and its subdifferential, ultimately forming a nonexpansive PDMM operator $\mathbf{T}_{P,\rho}$ to use as part of a classic fixed point iterative algorithm.

As the operator $\mathbf{I} - \mathbf{P}$ is continuous, $\ker(\mathbf{I} - \mathbf{P})$ is a closed subspace. It follows from (Example 1.25) [34] that $\iota_{\ker(\mathbf{I}-\mathbf{P})} \in \Gamma_0$. By combining this point with Theorem 13.37 and Proposition 13.3, both of [34], and the assumption that $f \in \Gamma_0$, it follows that $\boldsymbol{\lambda}^*$ is a minimiser of (5) if and only if (Theorem 16.3) [34]

$$\mathbf{0} \in \mathbf{C}\partial f^* \left(\mathbf{C}^T \boldsymbol{\lambda}^*\right) - \mathbf{d} + \partial \iota_{\ker(\mathbf{I}-\mathbf{P})}\left(\boldsymbol{\lambda}^*\right), \quad (6)$$

where $\partial \iota_{\ker(\mathbf{I}-\mathbf{P})}$ is a normal cone operator [39].

Given the maximal monotonicity of $\mathbf{T}_1 = \mathbf{C}\partial f^*\left(\mathbf{C}^T\right) - \mathbf{d}$ and of the normal cone operator $\mathbf{T}_2 = \partial \iota_{\ker(\mathbf{I}-\mathbf{P})}$, it follows that a zero-point of (6) can be found via a range of classical operator splitting methods (see [40] for an overview). In this particular instance, we consider the use of Peaceman-Rachford (PR) splitting to construct a nonexpansive PDMM operator which can be used to rephrase the zero-point condition in (6) as a more familiar fixed-point condition. A derivation of this equivalence, as demonstrated in (Section 7.3) [39], is given by

$$\mathbf{0} \in (\mathbf{I} + \rho \mathbf{T}_2)(\boldsymbol{\lambda}) - (\mathbf{I} - \rho \mathbf{T}_1)(\boldsymbol{\lambda})$$
$$\mathbf{0} \in (\mathbf{I} + \rho \mathbf{T}_2)(\boldsymbol{\lambda}) - \mathbf{R}_{\mathbf{T}_1,\rho} \circ (\mathbf{I} + \rho \mathbf{T}_1)(\boldsymbol{\lambda}) \quad (7a)$$
$$\mathbf{0} \in (\mathbf{I} + \rho \mathbf{T}_2) \circ \mathbf{J}_{\mathbf{T}_1,\rho}(\mathbf{z}) - \mathbf{R}_{\mathbf{T}_1,\rho}(\mathbf{z}), \; \boldsymbol{\lambda} = \mathbf{J}_{\mathbf{T}_1,\rho}(\mathbf{z})$$
$$\mathbf{J}_{\mathbf{T}_2,\rho} \circ \mathbf{R}_{\mathbf{T}_1,\rho}(\mathbf{z}) = \mathbf{J}_{\mathbf{T}_1,\rho}(\mathbf{z}), \; \boldsymbol{\lambda} = \mathbf{J}_{\mathbf{T}_1,\rho}(\mathbf{z}) \quad (7b)$$
$$\mathbf{R}_{\mathbf{T}_2,\rho} \circ \mathbf{R}_{\mathbf{T}_1,\rho}(\mathbf{z}) = \mathbf{z}, \; \boldsymbol{\lambda} = \mathbf{J}_{\mathbf{T}_1,\rho}(\mathbf{z}), \quad (7c)$$

where $\mathbf{R}_{\mathbf{T}_i,\rho}$ and $\mathbf{J}_{\mathbf{T}_i,\rho}$ are the reflected resolvent and resolvent operators of $\mathbf{T}_i$ respectively. Here, (7a) uses the maximality of $\mathbf{T}_1$ [39], (7b) uses that the reflected resolvent of a maximal monotone operator is single valued, and (7c) uses the definition of the reflected resolvent. The introduced *auxiliary* $\mathbf{z}$ variables will be referred to as such from here on out.

As both $\mathbf{T}_1$ and $\mathbf{T}_2$ are maximal monotone operators, $\mathbf{J}_{\mathbf{T}_1,\rho}$ and $\mathbf{J}_{\mathbf{T}_2,\rho}$ are both firmly nonexpansive. By Proposition 4.2 in [34] it follows that $\mathbf{R}_{\mathbf{T}_1,\rho}$ and $\mathbf{R}_{\mathbf{T}_2,\rho}$ are nonexpansive.

The resulting PDMM operator, given by

$$\mathbf{T}_{P,\rho} = \mathbf{R}_{\mathbf{T}_2,\rho} \circ \mathbf{R}_{\mathbf{T}_1,\rho}, \quad (8)$$

is nonexpansive and will be used repeatedly throughout this work. Importantly, the nonexpansive nature of $\mathbf{T}_{P,\rho}$ allows us to utilise fixed-point schemes to ultimately solve (1).

*C. On the link with the Primal Dual Method of Multipliers*

We will now demonstrate how PDMM, as defined in [31], can be linked with classical monotone operator splitting theory. For this purpose we will consider the fixed-point iteration of $\mathbf{T}_{P,\rho}$ given by

$$\mathbf{z}^{(k+1)} = \mathbf{T}_{P,\rho}\left(\mathbf{z}^{(k)}\right) = \mathbf{R}_{\mathbf{T}_2,\rho} \circ \mathbf{R}_{\mathbf{T}_1,\rho}\left(\mathbf{z}^{(k)}\right). \quad (9)$$

To aid in the aforementioned relationship, the evaluation of the reflected resolvent operators $\mathbf{R}_{\mathbf{T}_1,\rho}$ and $\mathbf{R}_{\mathbf{T}_2,\rho}$ are outlined in the following Lemmas.

**Lemma IV.1.** $\mathbf{y}^{(k+1)} = \mathbf{R}_{\mathbf{T}_1,\rho}\left(\mathbf{z}^{(k)}\right)$ *can be computed as*

$$\mathbf{x}^{(k+1)} = \arg\min_{\mathbf{x}} \left(f(\mathbf{x}) - \left\langle \mathbf{C}^T \mathbf{z}^{(k)}, \mathbf{x} \right\rangle + \frac{\rho}{2}\|\mathbf{C}\mathbf{x} - \mathbf{d}\|^2\right)$$
$$\boldsymbol{\lambda}^{(k+1)} = \mathbf{z}^{(k)} - \rho\left(\mathbf{C}\mathbf{x}^{(k+1)} - \mathbf{d}\right)$$
$$\mathbf{y}^{(k+1)} = 2\boldsymbol{\lambda}^{(k+1)} - \mathbf{z}^{(k)}$$

*Proof.* As $\mathbf{R}_{\mathbf{T}_1,\rho} = 2\mathbf{J}_{\mathbf{T}_1,\rho} - \mathbf{I}$, we begin by defining a method for computing the update

$$\boldsymbol{\lambda}^{(k+1)} = \mathbf{J}_{\mathbf{T}_1,\rho}\left(\mathbf{z}^{(k)}\right).$$

Firstly, by the definition of the resolvent of an operator,

$$\boldsymbol{\lambda}^{(k+1)} = (\mathbf{I} + \rho \mathbf{T}_1)^{-1}\left(\mathbf{z}^{(k)}\right)$$
$$\mathbf{z}^{(k)} \in (\mathbf{I} + \rho \mathbf{T}_1)\left(\boldsymbol{\lambda}^{(k+1)}\right)$$
$$\boldsymbol{\lambda}^{(k+1)} \in \mathbf{z}^{(k)} - \rho \mathbf{T}_1\left(\boldsymbol{\lambda}^{(k+1)}\right).$$

From the definition of the operator $\mathbf{T}_1$, it follows that

$$\boldsymbol{\lambda}^{(k+1)} \in \mathbf{z}^{(k)} - \rho\left(\mathbf{C}\partial f^*\left(\mathbf{C}^T \boldsymbol{\lambda}^{(k+1)}\right) - \mathbf{d}\right).$$

Let $\mathbf{x} \in \partial f^*\left(\mathbf{C}^T \boldsymbol{\lambda}\right)$. For $f \in \Gamma_0$, it follows from Proposition 16.10 [34], that $\mathbf{x} \in \partial f^*\left(\mathbf{C}^T \boldsymbol{\lambda}\right) \iff \partial f(\mathbf{x}) \ni \mathbf{C}^T \boldsymbol{\lambda}$ so that

$$\boldsymbol{\lambda}^{(k+1)} = \mathbf{z}^{(k)} - \rho\left(\mathbf{C}\mathbf{x}^{(k+1)} - \mathbf{d}\right)$$
$$\mathbf{C}^T \boldsymbol{\lambda}^{(k+1)} \in \partial f\left(\mathbf{x}^{(k+1)}\right). \quad (10)$$

Thus, $\mathbf{x}^{(k+1)}$ can be computed as

$$\mathbf{0} \in \partial f\left(\mathbf{x}^{(k+1)}\right) - \mathbf{C}^T\left(\mathbf{z}^{(k)} - \rho\left(\mathbf{C}\mathbf{x}^{(k+1)} - \mathbf{d}\right)\right) \quad (11)$$

or equivalently as

$$\mathbf{x}^{(k+1)} = \arg\min_{\mathbf{x}}\left(f(\mathbf{x}) - \left\langle \mathbf{C}^T \mathbf{z}^{(k)}, \mathbf{x} \right\rangle + \frac{\rho}{2}\|\mathbf{C}\mathbf{x} - \mathbf{d}\|^2\right) \quad (12)$$

Combining (10) and (12) with the fact that $\mathbf{y}^{(k+1)} = (2\mathbf{J}_{\mathbf{T}_1,\rho} - \mathbf{I})\left(\mathbf{z}^{(k)}\right) = 2\boldsymbol{\lambda}^{(k+1)} - \mathbf{z}^{(k)}$ completes the proof. $\square$

**Remark 2.** *Note that due to the block diagonal structure of $\mathbf{C}$ and the separability of $f$, the computation of this reflected resolvent can be performed in parallel across nodes.*

**Lemma IV.2.** $\mathbf{z}^{(k+1)} = \mathbf{R}_{\mathbf{T}_2,\rho}\left(\mathbf{y}^{(k+1)}\right)$ *can be computed as* $\mathbf{z}^{(k+1)} = \mathbf{P}\mathbf{y}^{(k+1)}$.

*Proof.* As in Lemma IV.1, we begin by utilising the definition of the reflected resolvent where $\mathbf{R}_{\mathbf{T}_2,\rho} = 2\mathbf{J}_{\mathbf{T}_2,\rho} - \mathbf{I}$. We therefore shift our attention to the computation of $\mathbf{J}_{\mathbf{T}_2,\rho}$.

The resolvent of a normal cone operator can be computed via

$$\mathbf{J}_{\mathbf{T}_2,\rho}\left(\mathbf{y}^{(k+1)}\right) = \arg\min_{\mathbf{u}}\left(\iota_{\ker(\mathbf{I}-\mathbf{P})}(\mathbf{u}) + \frac{1}{2\rho}\|\mathbf{u} - \mathbf{y}^{(k+1)}\|^2\right)$$
$$= \arg\min_{\mathbf{u}=\mathbf{P}\mathbf{u}}\left(\|\mathbf{u} - \mathbf{y}^{(k+1)}\|^2\right).$$

By inspection, the solution of this problem is given by the projection onto the set of feasible $\mathbf{u}$ vectors. As such,

$$\mathbf{J}_{\mathbf{T}_2,\rho}\left(\mathbf{y}^{(k+1)}\right) = \prod_{\ker(\mathbf{I}-\mathbf{P})} \mathbf{y}^{(k+1)}.$$

It follows that the reflected resolvent can be computed as

$$\mathbf{z}^{(k+1)} = \left(2 \prod_{\ker(\mathbf{I}-\mathbf{P})} - \mathbf{I}\right) \mathbf{y}^{(k+1)}$$
$$= \mathbf{P}\mathbf{y}^{(k+1)},$$

completing the proof. □

**Remark 3.** *Note that this permutation operation is equivalent to a sharing of data between neighbouring nodes and is therefore distributable.*

Utilising Lemmas IV.1 and IV.2 it follows that

$$\mathbf{T}_{P,\rho} = \mathbf{P} \circ \mathbf{R}_{\mathbf{T}_1,\rho},$$

so that (9) can be computed via Algorithm 1.

---

**Algorithm 1** Primal-Dual Method of Multipliers (PDMM)

1: **Initialise:** $\mathbf{z}^{(0)} \in \mathbb{R}^{M_E}$
2: **for** k=0,..., **do**
3: $\quad \mathbf{x}^{(k+1)} = \underset{\mathbf{x}}{\arg\min}\left(f(\mathbf{x}) - \langle \mathbf{C}^T \mathbf{z}^{(k)}, \mathbf{x}\rangle + \frac{\rho}{2}\|\mathbf{C}\mathbf{x}-\mathbf{d}\|^2\right)$
4: $\quad \boldsymbol{\lambda}^{(k+1)} = \mathbf{z}^{(k)} - \rho\left(\mathbf{C}\mathbf{x}^{(k+1)} - \mathbf{d}\right)$
5: $\quad \mathbf{y}^{(k+1)} = 2\boldsymbol{\lambda}^{(k+1)} - \mathbf{z}^{(k)}$
6: $\quad \mathbf{z}^{(k+1)} = \mathbf{P}\mathbf{y}^{(k+1)}$
7: **end for**

---

By noting that $\mathbf{z}^{(k+1)} = \mathbf{P}\left(\boldsymbol{\lambda}^{(k+1)} - \rho\left(\mathbf{C}\mathbf{x}^{(k+1)} - \mathbf{d}\right)\right)$, the dependence on $\mathbf{y}^{(k+1)}$ and $\mathbf{z}^{(k+1)}$ can be removed, reducing the scheme to that given in Algorithm 2. This algorithm

---

**Algorithm 2** Simplified PDMM

1: **Initialise:** $\boldsymbol{\lambda}^{(0)} \in \mathbb{R}^{M_E}$, $\mathbf{x}^{(0)} \in \mathbb{R}^{M_V}$
2: **for** k=0,..., **do**
3: $\quad \mathbf{x}^{(k+1)} = \underset{\mathbf{x}}{\arg\min}\big(f(\mathbf{x}) - \langle \mathbf{C}^T \mathbf{P}\boldsymbol{\lambda}^{(k)}, \mathbf{x}\rangle +$
$\quad\quad \frac{\rho}{2}\|\mathbf{C}\mathbf{x} + \mathbf{P}\mathbf{C}\mathbf{x}^{(k)} - 2\mathbf{d}\|^2\big)$
4: $\quad \boldsymbol{\lambda}^{(k+1)} = \mathbf{P}\boldsymbol{\lambda}^{(k)} - \rho\left(\mathbf{C}\mathbf{x}^{(k+1)} + \mathbf{P}\mathbf{C}\mathbf{x}^{(k)} - 2\mathbf{d}\right)$
5: **end for**

---

is identical to a particular instance of PDMM proposed in [31]. Thus, PDMM is equivalent to the fixed-point iteration of the PR splitting of the extended dual problem, linking the approach with a plethora of existing algorithms within the literature [21], [23], [24], [26].

The distributed nature of PDMM can be more easily visualised in Algorithm 3 where by we have utilised the definitions of $\mathbf{C}$ and $\mathbf{d}$. Here the notation $\mathbf{Node}_j \leftarrow \mathbf{Node}_i(\bullet)$ indicates the transmission of data from node $i$ to node $j$.

**Remark 4.** *The form of the distributed PDMM algorithm lends itself to implementation in practical networks. Each iteration of the algorithm only requires one-way transmission between neighbouring nodes. Thus, no direct collaboration is required during the computations of each iteration.*

---

**Algorithm 3** Distributed PDMM

1: **Initialise:** $\mathbf{z}^{(0)} \in \mathbb{R}^{M_E}$
2: **for** k=0,..., **do**
3: $\quad$ **for all** $i \in V$ **do** ▷ Primal Update
4: $\quad\quad \mathbf{x}_i^{(k+1)} = \arg\min_{\mathbf{x}_i}\Big(f_i(\mathbf{x}_i) +$
$\quad\quad \sum_{j \in \mathcal{N}(i)}\left(-\left\langle \mathbf{A}_{i|j}^T \mathbf{z}_{i|j}^{(k)}, \mathbf{x}_i\right\rangle + \frac{\rho}{2}\|\mathbf{A}_{i|j}\mathbf{x}_i - \frac{\mathbf{b}_{ij}}{2}\|^2\right)\Big)$
5: $\quad\quad$ **for all** $j \in \mathcal{N}(i)$ **do** ▷ Dual Update
6: $\quad\quad\quad \boldsymbol{\lambda}_{i|j}^{(k+1)} = \mathbf{z}_{i|j}^{(k)} - \rho\left(\mathbf{A}_{i|j}\mathbf{x}_i^{(k+1)} - \frac{\mathbf{b}_{ij}}{2}\right)$
7: $\quad\quad\quad \mathbf{z}_{i|j}^{(k+\frac{1}{2})} = 2\boldsymbol{\lambda}_{i|j}^{(k+1)} - \mathbf{z}_{i|j}^{(k)}$
8: $\quad\quad$ **end for**
9: $\quad$ **end for**
10: $\quad$ **for all** $i \in V, j \in \mathcal{N}(i)$ **do** ▷ Transmit Variables
11: $\quad\quad \mathbf{Node}_j \leftarrow \mathbf{Node}_i(\mathbf{z}_{i|j}^{(k+\frac{1}{2})})$
12: $\quad$ **end for**
13: $\quad$ **for all** $i \in V, j \in \mathcal{N}(i)$ **do** ▷ Auxiliary Update
14: $\quad\quad \mathbf{z}_{i|j}^{(k+1)} = \mathbf{z}_{j|i}^{(k+\frac{1}{2})}$
15: $\quad$ **end for**
16: **end for**

---

*D. On the Link with the Distributed Alternating Direction Method of Multipliers*

While in [31] it was suggested that PDMM and the alternating direction method of multipliers (ADMM) were fundamentally different due to their differences in derivation, in the following we demonstrate how they are more closely related than first thought. This link is masked via the change of variables typically used in the updating scheme for ADMM and PDMM (see [23] Sec. 3 and [31] Sec. 4 respectively for such representations). For this purpose we re-derive an ADMM variant from the perspective of monotone operator theory.

Consider again the problem given in (1). By introducing the directed edge variables $\mathbf{y}_{i|j}, \mathbf{y}_{j|i} \in \mathbb{R}^{M_{ij}}$ $\forall (i,j) \in E$, an equivalent version of (1) is given by

$$\min_{\mathbf{x}} \sum_{i \in V} f_i(\mathbf{x_i})$$
$$\text{s.t} \quad \left.\begin{array}{l}\mathbf{A}_{i|j}\mathbf{x}_i - \frac{\mathbf{b}_{ij}}{2} = \mathbf{y}_{i|j} \\ \mathbf{A}_{j|i}\mathbf{x}_j - \frac{\mathbf{b}_{ij}}{2} = \mathbf{y}_{j|i} \\ \mathbf{y}_{i|j} + \mathbf{y}_{i|j} = \mathbf{0}\end{array}\right\} \forall (i,j) \in E. \quad (13)$$

Defining the stacked vector $\mathbf{y} \in \mathbb{R}^{M_E}$ and adopting the matrices $\mathbf{C}$, $\mathbf{P}$ and $\mathbf{d}$ as per the derivation in Sec. IV-A, this can be more simply written as

$$\begin{aligned}\min_{\mathbf{x}} \quad & f(\mathbf{x}) + \iota_{\ker(\mathbf{I}+\mathbf{P})}(\mathbf{y}) \\ \text{s.t} \quad & \mathbf{C}\mathbf{x} - \mathbf{d} = \mathbf{y}.\end{aligned} \quad (14)$$

Here, the indicator function is used to capture the final set of equality constraints in (13). This is exactly in the form of a standard ADMM problem with the objective being the summation for two CCP functions of different variables with affine constraints coupling them together. The associated dual problem of (14) is given by

$$\min_{\boldsymbol{\lambda}_i \forall i \in V} f^*\left(\mathbf{C}^T \boldsymbol{\lambda}\right) - \mathbf{d}^T \boldsymbol{\lambda} + \iota^*_{\ker(\mathbf{I}+\mathbf{P})}(\boldsymbol{\lambda}), \quad (15)$$

where $\boldsymbol{\lambda}$, as in the case of PDMM, denotes the stacked vector of dual variables associated with the directed edges.

Comparing (15) and (6), we can note that the apparent difference in the dual problems is simply the use of either an indicator function, in the case of PDMM, or the conjugate of an indicator function in the case of ADMM. In actual fact these two functions are equivalent. This can be trivially seen by considering the definition of the Fenchel conjugate of an indicator function given by

$$\iota^*_{\ker(\mathbf{I}+\mathbf{P})}(\boldsymbol{\lambda}) = \sup_{\mathbf{y}}\left(\langle \mathbf{y},\boldsymbol{\lambda}\rangle - \iota_{\ker(\mathbf{I}+\mathbf{P})}\right)$$
$$= \begin{cases} 0 & \boldsymbol{\lambda} \in \mathrm{ran}\,(\mathbf{I}+\mathbf{P}) \\ -\infty & \text{otherwise.} \end{cases}$$

As $\mathrm{ran}\,(\mathbf{I}+\mathbf{P}) = \ker(\mathbf{I}-\mathbf{P})$, it follows that $\iota^*_{\ker(\mathbf{I}+\mathbf{P})} = \iota_{\ker(\mathbf{I}-\mathbf{P})}$. The problems in (5) and (15) are therefore identical.

ADMM is equivalent to applying Douglas Rachford (DR) splitting [22] to (15). As DR splitting is equivalent to a half averaged form of PR splitting [34], the operator form of ADMM is given by $\mathbf{T}_{A,\rho} = \frac{1}{2}(\mathbf{I}+\mathbf{T}_{P,\rho})$. In this manner, despite their differences in earlier derivations, the two approaches are fundamentally linked. Within the literature, PDMM could therefore also be referred to as a particular instance of *generalised* [41] or *relaxed* ADMM [42].

## V. GENERAL CONVERGENCE RESULTS FOR PDMM

Having linked PDMM with more classical monotone based approaches, we now move to demonstrate a set of convergence results of the algorithm under different assumptions on the underlying optimisation problem. These results include a sufficient condition for strong primal convergence and a reiteration of the asymptotic rate bound for averaged PDMM variants, which provides guaranteed convergence for all $f \in \Gamma_0$.

### A. Sufficient Conditions for the Optimality of the PDMM Variables in the Limit

Firstly we show how either the differentiability or the strict convexity of $f$ is sufficient for the optimality of the limit points of the primal or dual variables of PDMM respectively. Strengthening this, when $f$ is both strictly convex and differentiable we will demonstrate how PDMM guarantees convergence of the auxiliary variables to a particular set, in turn guaranteeing primal convergence.

Consider the PDMM operator $\mathbf{T}_{P,\rho}$ given in (8). Assume that $\mathrm{fix}\,(\mathbf{T}_{P,\rho}) \neq \emptyset$ and define a point $\mathbf{z}^* \in \mathrm{fix}\,(\mathbf{T}_{P,\rho})$. We will denote by $\mathbf{x}^*$ and $\boldsymbol{\lambda}^*$ the corresponding primal and dual vectors generated by the PDMM algorithm corresponding to $\mathbf{z}^*$. Consider the sequence $(\mathbf{z}^{(k)})_{k \in \mathbb{Z}_+}$ as defined by (9). It follows that, for a given iteration $k \in \mathbb{Z}_+$,

$$\|\mathbf{z}^{(k)} - \mathbf{z}^*\|^2 = \|\mathbf{P}(\mathbf{y}^{(k)} - \mathbf{y}^*)\|^2$$
$$= \|\mathbf{y}^{(k)} - \mathbf{y}^*\|^2$$
$$= \|\mathbf{z}^{(k-1)} - \mathbf{z}^*\|^2 - 4\left\langle \boldsymbol{\lambda}^{(k)} - \boldsymbol{\lambda}^*, \rho \mathbf{C}\left(\mathbf{x}^{(k)} - \mathbf{x}^*\right)\right\rangle,$$

where the final line uses the polarization identity of norms (Lemma 2.12) [34] and the definition of the dual update step in Lemma IV.1. As $\mathbf{x} \in \partial f^*(\mathbf{C}^T\boldsymbol{\lambda})$, by Proposition 16.10 [34], $\mathbf{C}^T\boldsymbol{\lambda} \in \partial f(\mathbf{x})$. It follows that

$$\left\langle \boldsymbol{\lambda}^{(k)} - \boldsymbol{\lambda}^*, \mathbf{C}\left(\mathbf{x}^{(k)} - \mathbf{x}^*\right)\right\rangle =$$
$$\left\langle \partial f\left(\mathbf{x}^{(k)}\right) - \partial f\left(\mathbf{x}^*\right), \mathbf{x}^{(k)} - \mathbf{x}^*\right\rangle =$$
$$\left\langle \boldsymbol{\lambda}^{(k)} - \boldsymbol{\lambda}^*, \mathbf{C}\partial f^*\left(\mathbf{C}^T\boldsymbol{\lambda}^{(k)}\right) - \mathbf{C}\partial f^*\left(\mathbf{C}^T\boldsymbol{\lambda}^*\right)\right\rangle.$$

If $\partial f$ is strictly monotone

$$\left\langle \partial f\left(\mathbf{x}^{(k)}\right) - \partial f\left(\mathbf{x}^*\right), \mathbf{x}^{(k)} - \mathbf{x}^*\right\rangle > 0 \;\forall \mathbf{x}^{(k)} \neq \mathbf{x}^*. \quad (16)$$

Similarly, if $\partial f^*$ is strictly monotone

$$\left\langle \boldsymbol{\lambda}^{(k)} - \boldsymbol{\lambda}^*, \mathbf{C}f^*\left(\mathbf{C}^T\boldsymbol{\lambda}^{(k)}\right) - \mathbf{C}f^*\left(\mathbf{C}^T\boldsymbol{\lambda}^*\right)\right\rangle > 0$$
$$\forall \Pi_{\mathrm{ran}(\mathbf{C})}\left(\boldsymbol{\lambda}^{(k)}\right) \neq \Pi_{\mathrm{ran}(\mathbf{C})}(\boldsymbol{\lambda}^*). \quad (17)$$

From Rockafellar ([43] Theorem 23.6), a function $f \in \Gamma_0$ is strictly convex if and only if $f^*$ is differentiable. Given that $f^{**} = f$ for $f \in \Gamma_0$, the strict monotonicity of $\partial f^*$ holds if and only if $f$ is differentiable. If $f$ is differentiable or strictly convex, we now show how the limits of the primal and dual variables are optimal respectively.

**Lemma V.1.** *In the limit if*
*(i) $f$ is strictly convex then $\lim_{k\to\infty}\mathbf{x}^{(k)} = \mathbf{x}^*$*
*(ii) $\partial f = \nabla f$ then $\lim_{k\to\infty} \Pi_{\mathrm{ran}(\mathbf{C})}\left(\boldsymbol{\lambda}^{(k)}\right) = \Pi_{\mathrm{ran}(\mathbf{C})}(\boldsymbol{\lambda}^*)$*

*Proof.* Consider the first case *(i)*. Due to the strict convexity of $f$ and the convexity of the domain of (1), $\mathbf{x}^* \in \mathbf{X}^*$ is unique where $\mathbf{X}^*$ denotes the set of primal optimal vectors.

Under such an assumption, from (16), it follows that the sequence $(\|\mathbf{z}^{(k)} - \mathbf{z}^*\|^2)_{k \in \mathbb{Z}_+}$ is strictly nonexpansive unless $\mathbf{x}^{(k)} = \mathbf{x}^*$. Thus $\lim_{k\to\infty}\mathbf{x}^{(k)} = \mathbf{x}^*$ in this case.

For the second case *(ii)* we have a similar result. Denote by $\boldsymbol{\Lambda}^*$ the set of optimal dual variables. Due to the differentiability of $f$, $\mathbf{C}\partial f^*(\mathbf{C}^T\bullet)$ is strictly monotone over $\mathrm{ran}\,(\mathbf{C})$. It follows that $\forall \boldsymbol{\lambda} \in \boldsymbol{\Lambda}^*$, $\Pi_{\mathrm{ran}(\mathbf{C})}(\boldsymbol{\lambda}) = \Pi_{\mathrm{ran}(\mathbf{C})}(\boldsymbol{\lambda}^*)$ is unique.

From (17), the sequence $(\|\mathbf{z}^{(k)} - \mathbf{z}^*\|^2)_{k \in \mathbb{Z}_+}$ is strictly nonexpansive unless $\Pi_{\mathrm{ran}(\mathbf{C})}\left(\boldsymbol{\lambda}^{(k)}\right) = \Pi_{\mathrm{ran}(\mathbf{C})}(\boldsymbol{\lambda}^*)$. Thus $\lim_{k\to\infty}\Pi_{\mathrm{ran}(\mathbf{C})}\left(\boldsymbol{\lambda}^{(k)}\right) = \Pi_{\mathrm{ran}(\mathbf{C})}(\boldsymbol{\lambda}^*)$ in this case. $\square$

The importance of Lemma V.1 is that unlike the original derivation of PDMM [31], here we postulate that primal convergence cannot be guaranteed for all $f \in \Gamma_0$ without additional modifications to the algorithm. This modification may take the form of operator averaging, as introduced in Section V-C although other options may also be sufficient. This postulation is verified in Sec. VII-B by demonstrating an instance where $f \in \Gamma_0$ is insufficient for primal optimality.

### B. Sufficient Conditions for Auxiliary Convergence of PDMM

We will now demonstrate that if the function $f$ is both differentiable and strictly convex we can strengthen the result of Sec. V-A to guarantee auxiliary convergence. This also guarantees primal convergence for the same function class.

Consider the one and two-step PDMM updates given by

$$\begin{aligned}\mathbf{z}^{(k+1)} &= \mathbf{P}\left(\mathbf{z}^{(k)} - 2\rho\left(\mathbf{C}\mathbf{x}^{(k+1)} - \mathbf{d}\right)\right) \\ \mathbf{z}^{(k+2)} &= \mathbf{P}\left(\mathbf{z}^{(k+1)} - 2\rho\left(\mathbf{C}\mathbf{x}^{(k+2)} - \mathbf{d}\right)\right) \\ &= \mathbf{z}^{(k)} - 2\rho\left(\mathbf{P}\mathbf{C}\mathbf{x}^{(k+2)} + \mathbf{C}\mathbf{x}^{(k+1)} - 2\mathbf{d}\right).\end{aligned} \quad (18)$$

From (18) it is clear that for $k$ even and for any given $\mathbf{z}^0$,

$$\prod_{\ker(\mathbf{C})\cap\ker(\mathbf{P}\mathbf{C})} \mathbf{z}^{(k+2)} = \prod_{\ker(\mathbf{C})\cap\ker(\mathbf{P}\mathbf{C})} \mathbf{z}^{(0)}.$$

Similarly for $k$ odd

$$\prod_{\ker(\mathbf{C})\cap\ker(\mathbf{P}\mathbf{C})} \mathbf{z}^{(k+2)} = \prod_{\ker(\mathbf{C})\cap\ker(\mathbf{P}\mathbf{C})} \mathbf{T}_{P,\rho}\left(\mathbf{z}^{(0)}\right).$$

The auxiliary variables therefore contain a nonconvergent component determined by the vector $\mathbf{z}_0$. We therefore define a measure of auxiliary error, given by $\|\mathbf{z}^{(k+2)} - \tilde{\mathbf{z}}_0\|^2$ where the vector $\tilde{\mathbf{z}}_0$ is defined to remove the dependence on this nonconvergent component.

Our particular choice of $\tilde{\mathbf{z}}_0$ is derived as follows. Assume that $f$ is both differentiable and strictly convex. Consider the two-step PDMM operator $\mathbf{T}_{P,\rho}^2 = \mathbf{T}_{P,\rho} \circ \mathbf{T}_{P,\rho}$, and note that $\mathbf{z}^* \in \text{fix}(\mathbf{T}_{P,\rho}) \subseteq \text{fix}(\mathbf{T}_{P,\rho}^2)$. From (18) and as such a point produces a unique primal optimal update as $f$ is strictly convex, it follows that

$$\begin{aligned}\mathbf{z}^* &= \mathbf{P}\left(\mathbf{z}^* - 2\rho\left(\mathbf{C}\mathbf{x}^* - \mathbf{d}\right)\right) \\ &= \mathbf{z}^* - 2\rho\left(\mathbf{P}\mathbf{C}\mathbf{x}^* + \mathbf{C}\mathbf{x}^* - 2\mathbf{d}\right) = \mathbf{z}^*,\end{aligned}$$

where $\mathbf{P}\mathbf{C}\mathbf{x}^* + \mathbf{C}\mathbf{x}^* - 2\mathbf{d} = \mathbf{0}$ so that $\mathbf{x}^*$ is primal feasible.

Additionally, from the primal updates of PDMM and the differentiability of $f$, for a single PDMM iteration

$$\mathbf{x}^* = \nabla f(\mathbf{x}^*) - \mathbf{C}^T\mathbf{z}^* + \rho\mathbf{C}^T\left(\mathbf{C}\mathbf{x}^* - \mathbf{d}\right).$$

while for a second iteration

$$\begin{aligned}\mathbf{x}^* &= \nabla f(\mathbf{x}^*) - \mathbf{C}^T\mathbf{z}^* + \rho\mathbf{C}^T\left(\mathbf{C}\mathbf{x}^* - \mathbf{d}\right) \\ &= \nabla f(\mathbf{x}^*) - \mathbf{C}^T\mathbf{P}\mathbf{z}^* + \rho\mathbf{C}^T\mathbf{P}\left(\mathbf{C}\mathbf{x}^* - \mathbf{d}\right),\end{aligned}$$

where again we have used the fact that $\mathbf{P}\mathbf{C}\mathbf{x}^* + \mathbf{C}\mathbf{x}^* - 2\mathbf{d} = \mathbf{0}$. By stacking these two equalities, it follows that

$$\begin{bmatrix}\mathbf{C} & \mathbf{P}\mathbf{C}\end{bmatrix}^T \mathbf{z}^* = \begin{bmatrix}\nabla f(\mathbf{x}^*) - \mathbf{x}^* + \rho\mathbf{C}^T\left(\mathbf{C}\mathbf{x}^* - \mathbf{d}\right) \\ \nabla f(\mathbf{x}^*) - \mathbf{x}^* + \rho\mathbf{C}^T\mathbf{P}\left(\mathbf{C}\mathbf{x}^* - \mathbf{d}\right)\end{bmatrix} = \hat{\mathbf{g}},$$

and thus

$$\prod_{\text{ran}(\mathbf{C})+\text{ran}(\mathbf{P}\mathbf{C})} \mathbf{z}^* = \left(\begin{bmatrix}\mathbf{C} & \mathbf{P}\mathbf{C}\end{bmatrix}^T\right)^\dagger \hat{\mathbf{g}},$$

where $\bullet^\dagger$ denotes the Moore-Penrose pseudo inverse.

We can therefore define the unique vector

$$\tilde{\mathbf{z}}_0 = \left(\begin{bmatrix}\mathbf{C} & \mathbf{P}\mathbf{C}\end{bmatrix}^T\right)^\dagger \hat{\mathbf{g}} + \prod_{\ker(\mathbf{C})\cap\ker(\mathbf{P}\mathbf{C})} \mathbf{z}^{(0)}. \quad (19)$$

By again considering the PDMM two-step operator $\mathbf{T}_{P,\rho}^2$, it follows from nonexpansiveness that, for $k$ even,

$$\begin{aligned}\|\mathbf{z}^{(k+2)} - \tilde{\mathbf{z}}_0\|^2 &= \|\prod_{\text{ran}(\mathbf{C})+\text{ran}(\mathbf{P}\mathbf{C})}\left(\mathbf{z}^{(k+2)} - \tilde{\mathbf{z}}_0\right)\|^2 \\ &\leq \|\prod_{\text{ran}(\mathbf{C})+\text{ran}(\mathbf{P}\mathbf{C})}\left(\mathbf{z}^{(k)} - \tilde{\mathbf{z}}_0\right)\|^2 \\ &= \|\mathbf{z}^{(k)} - \tilde{\mathbf{z}}_0\|^2.\end{aligned} \quad (20)$$

Thus the sequence $\left(\mathbf{z}^{(2k)}\right)_{k\in\mathbb{Z}_+}$ is Fejér monotone with respect to the point $\tilde{\mathbf{z}}_0$. It follows from (Proposition 5.4) [34] that the sequence of auxiliary errors is both bounded and converges. Note that for odd iterations the vector $\tilde{\mathbf{z}}_1 = \mathbf{T}_{P,\rho}(\tilde{\mathbf{z}}_0)$ can be used to define a similar error.

Using the proposed auxiliary error, the limit points of the sequence of two-step auxiliary variables can be derived. We do so by noting that for differentiable and strictly convex $f$, Lemma V.1 implies that

$$\begin{aligned}\lim_{k\to\infty} \mathbf{x}^{(2k)} &= \lim_{k\to\infty} \mathbf{x}^{(2k+1)} = \mathbf{x}^* \\ \lim_{k\to\infty} \prod_{\text{ran}(\mathbf{C})}\left(\boldsymbol{\lambda}^{(2k)}\right) &= \prod_{\text{ran}(\mathbf{C})}\left(\tilde{\boldsymbol{\lambda}}_0\right) \\ \lim_{k\to\infty} \prod_{\text{ran}(\mathbf{C})}\left(\boldsymbol{\lambda}^{(2k+1)}\right) &= \prod_{\text{ran}(\mathbf{C})}\left(\tilde{\boldsymbol{\lambda}}_1\right),\end{aligned} \quad (21)$$

where $\tilde{\boldsymbol{\lambda}}_0 = \mathbf{J}_{\mathbf{T}_1,\rho}(\tilde{\mathbf{z}}_0)$ and $\tilde{\boldsymbol{\lambda}}_1 = \mathbf{J}_{\mathbf{T}_1,\rho} \circ \mathbf{T}_{P,\rho}(\tilde{\mathbf{z}}_0)$. From the definition of the reflected resolvent it follows that

$$\begin{aligned}\lim_{k\to\infty} \prod_{\text{ran}(\mathbf{C})}\left(\mathbf{z}^{(2k)}\right) &= \prod_{\text{ran}(\mathbf{C})}\left(\tilde{\mathbf{z}}_0\right) \\ \lim_{k\to\infty} \prod_{\text{ran}(\mathbf{C})}\left(\mathbf{z}^{(2k+1)}\right) &= \prod_{\text{ran}(\mathbf{C})}\left(\tilde{\mathbf{z}}_1\right),\end{aligned} \quad (22)$$

where $\tilde{\mathbf{z}}_1 = \mathbf{T}_{P,\rho}(\tilde{\mathbf{z}}_0)$. Combining (22) with (21) and the one-step PDMM update given in (18), it follows that

$$\begin{aligned}\mathbf{0} &= \lim_{k\to\infty} \prod_{\text{ran}(\mathbf{C})}\left(\mathbf{z}^{(2k+1)} - \tilde{\mathbf{z}}_1\right) \\ &= \lim_{k\to\infty} \prod_{\text{ran}(\mathbf{C})} \mathbf{P}\left(\mathbf{z}^{(2k)} - \tilde{\mathbf{z}}_0 - 2\rho\mathbf{C}\left(\mathbf{x}^{(k+1)} - \mathbf{x}^*\right)\right) \\ &= \lim_{k\to\infty} \mathbf{P} \prod_{\text{ran}(\mathbf{C})} \mathbf{P}\left(\mathbf{z}^{(2k)} - \tilde{\mathbf{z}}_0\right) \\ &= \lim_{k\to\infty} \prod_{\text{ran}(\mathbf{P}\mathbf{C})}\left(\mathbf{z}^{(2k)} - \tilde{\mathbf{z}}_0\right),\end{aligned} \quad (23)$$

where the second line uses Lemma IV.1 and the third line uses the fact that $\lim_{k\to\infty} \mathbf{x}^{(2k+1)} = \mathbf{x}^*$ and that $\mathbf{P}$ is full rank.

By then combining (22) and (23), it follows that

$$\lim_{k\to\infty} \prod_{\text{ran}(\mathbf{C})+\text{ran}(\mathbf{P}\mathbf{C})}\left(\mathbf{z}^{(2k)}\right) = \prod_{\text{ran}(\mathbf{C})+\text{ran}(\mathbf{P}\mathbf{C})}\left(\tilde{\mathbf{z}}_0\right),$$

so that $\lim_{k\to\infty} \|\mathbf{z}^{(2k)} - \tilde{\mathbf{z}}_0\|^2 = 0$. The vector $\tilde{\mathbf{z}}_0$ is thus the unique weak sequential cluster point of the sequence $\left(\mathbf{z}^{(2k)}\right)_{k\in\mathbb{Z}_+}$. Due to the finite dimensionality of the domain of $\mathbf{T}_{P,\rho}$, combining Theorem 5.5 and Lemma 2.51, both from [34], verifies that $\mathbf{z}^{(2k)} \to \tilde{\mathbf{z}}_0$ and $\mathbf{z}^{(2k+1)} \to \tilde{\mathbf{z}}_1$ where $\bullet \to \bullet$ is used to denote strong convergence.

Strong convergence of the primal variables follows by reconsidering the primal update equation given in Lem. IV.1. Specifically, by rearranging (11), it can be shown that

$$\mathbf{x}^{(k+1)} = \left(\nabla f + \rho\mathbf{C}^T\mathbf{C}\right)^{-1} \mathbf{C}^T\left(\mathbf{z}^{(k)} + \rho\mathbf{d}\right).$$

From the monotonicity of $\nabla f$, $\forall k \in \mathbb{Z}_+$ even,

$$\begin{aligned}\|\mathbf{x}^{(k+1)} - \mathbf{x}^*\|^2 &\leq \frac{1}{\rho\sigma_{\min\neq 0}^2(\mathbf{C})}\|\mathbf{C}^T\left(\mathbf{z}^{(k)} - \tilde{\mathbf{z}}_0\right)\|^2 \\ &\leq \frac{\sigma_{\max}^2(\mathbf{C})}{\rho\sigma_{\min\neq 0}^2(\mathbf{C})}\|\mathbf{z}^{(k)} - \tilde{\mathbf{z}}_0\|^2\end{aligned} \quad (24)$$

Thus, $\|\mathbf{x}^{(k+1)} - \mathbf{x}^*\|^2 \to 0$, $\mathbf{x}^{(k+1)} \to \mathbf{x}^*$ and the combination of differentiability and strict convexity is sufficient for strong primal convergence.

While this proof is sufficient for primal convergence, no information is provided about the rate at which this convergence occurs. However, as we will demonstrate in Secs. V-C and VI-A, with either the use of operator averaging or under stronger functional assumptions, such rates can be defined.

### C. Averaged PDMM Convergence

As with other operator splitting methods, PDMM can be combined with an averaging stage to guarantee convergence $\forall f \in \Gamma_0$, even those which do not satisfy the strict convexity or differentiability assumptions introduced in Sec. V-A. The general form of the averaged PDMM operator is given by

$$\mathbf{T}_{P,\rho,\alpha} = (1-\alpha)\mathbf{I} + \alpha \mathbf{T}_{P,\rho},$$

where the scalar $\alpha \in (0,1)$. In the particular case that $\alpha = \frac{1}{2}$, averaged PDMM is equivalent to ADMM, as was previously noted in Sec. IV-D. In this case, by Proposition 4.4 [34], the operator $\mathbf{T}_{P,\rho,\alpha}$ is firmly nonexpansive.

The fixed-point iteration of $\mathbf{T}_{P,\rho,\alpha}$ is therefore given by

$$\mathbf{z}^{(k+1)} = (1-\alpha)\mathbf{z}^{(k)} + \alpha \mathbf{T}_{P,\rho}\mathbf{z}^{(k)}.$$

This can be equivalently referred to as the $\alpha$-Krasnosel'skiĭ-Mann iteration [34] of the operator $\mathbf{T}_{P,\rho}$ which is a well documented method of guaranteeing convergence for nonexpansive operators. To demonstrate this, the following is a repetition of the guaranteed convergence of Krasnosel'skiĭ-Mann iterations, as demonstrated in Theorem 5.15 [34].

Consider again the sequence $(\mathbf{z}^{(k)})_{k \in \mathbb{Z}_+}$ and a vector $\mathbf{z}^* \in \mathrm{fix}(\mathbf{T}_{P,\rho})$. By the definition of the averaged PDMM operator, it follows that

$$\begin{aligned}\|\mathbf{z}^{(k+1)} - \mathbf{z}^*\|^2 &= \|\left((1-\alpha)\mathbf{I} + \alpha \mathbf{T}_{P,\rho}\right)\left(\mathbf{z}^{(k)}\right) - \mathbf{z}^*\|^2 \\ &= (1-\alpha)\|\mathbf{z}^{(k)} - \mathbf{z}^*\|^2 + \alpha\|\mathbf{T}_{P,\rho}\left(\mathbf{z}^{(k)}\right) - \mathbf{z}^*\|^2 \\ &\quad - \alpha(1-\alpha)\|\left(\mathbf{T}_{P,\rho} - \mathbf{I}\right)\left(\mathbf{z}^{(k)}\right)\|^2 \\ &\leq \|\mathbf{z}^{(k)} - \mathbf{z}^*\|^2 - \alpha(1-\alpha)\|\left(\mathbf{T}_{P,\rho} - \mathbf{I}\right)\left(\mathbf{z}^{(k)}\right)\|^2,\end{aligned} \quad (25)$$

where the first equality stems from Corollary 2.15 [34] and the final line uses the nonexpansiveness of $\mathbf{T}_{P,\rho}$. Recursively applying (25), if follows that

$$\sum_{i=1}^{k} \alpha(1-\alpha)\|\left(\mathbf{T}_{P,\rho} - \mathbf{I}\right)\left(\mathbf{z}^{(i)}\right)\|^2 \leq \|\mathbf{z}^0 - \mathbf{z}^*\|^2. \quad (26)$$

As shown in Theorem 5.15 [34],

$$\|\left(\mathbf{T}_{P,\rho} - \mathbf{I}\right)\left(\mathbf{z}^{(k)}\right)\|^2 \leq \|\left(\mathbf{T}_{P,\rho} - \mathbf{I}\right)\left(\mathbf{z}^{(k-1)}\right)\|^2,$$

so that the inequality in (26) can be lower bounded as

$$k\alpha(1-\alpha)\|\left(\mathbf{T}_{P,\rho} - \mathbf{I}\right)\left(\mathbf{z}^{(k)}\right)\|^2 \leq \|\mathbf{z}^0 - \mathbf{z}^*\|^2.$$

This leads to the final inequality,

$$\|\left(\mathbf{T}_{P,\rho} - \mathbf{I}\right)\left(\mathbf{z}^{(k)}\right)\|^2 \leq \frac{1}{k\alpha(1-\alpha)}\|\mathbf{z}^0 - \mathbf{z}^*\|^2.$$

It follows that the fixed-point residual $(\mathbf{T}_{P,\rho} - \mathbf{I})(\mathbf{z}^{(k)})$ converges at an asymptotic rate of $\mathcal{O}\left(\frac{1}{k}\right)$ and thus that $\mathbf{z}^{(k)}$ converges to $\mathrm{fix}(\mathbf{T}_{P,\rho})$. As such, operator averaging is a form of algorithm modification, alluded to in Sec. V-A, which guarantees primal convergence of PDMM for general $f \in \Gamma_0$.

## VI. GEOMETRIC CONVERGENCE AND DISTRIBUTED PARAMETER SELECTION

Under stronger functional assumptions, namely strong convexity and smoothness, we can demonstrate that a geometrically contracting upper bound exists for the primal error of PDMM. To complement this, we also introduce a finite-time distributed algorithm to select a step size to optimise this bound and improve the worst case convergence rate of PDMM.

### A. A Primal Geometric Convergence Bound for Strongly Convex and Smooth Functions

In the following we demonstrate that for strongly convex and smooth functions the primal variables of PDMM converge at a geometric rate. More formally we show that $\exists \, \epsilon \geq 0, \gamma \in [0, 1)$ so that

$$\forall k \in \mathbb{Z}_+, \, \|\mathbf{x}^{(k)} - \mathbf{x}^*\|^2 \leq \gamma^k \epsilon,$$

This process is broken down into two stages. Firstly we demonstrate how, under the assumption of strong convexity and smoothness, PDMM is contractive over a certain subspace. We then show how, for such "partially contractive" operators, a global convergence bound can be found by linking PDMM with the generalised alternating method of projections (GAP) [44] allowing us to derive the aforementioned $\gamma$ and $\epsilon$.

### B. Partially Contractive Nature of PDMM Over a Subspace

We begin by demonstrating the contractive nature of the PDMM operator over a subspace determined by the constraint matrix $\mathbf{C}$ by showing that the operator $\mathbf{T}_1$ is both strongly convex and Lipschitz continuous over this subspace. This is summarised in Lemma VI.1.

**Lemma VI.1.** *If $f$ is both $\mu$-strongly convex and $\beta$-smooth then $\mathbf{C}\partial f^*(\mathbf{C}^T \bullet)$ is*

(i) $\frac{\sigma_{\max}^2(\mathbf{C})}{\mu}$-*Lipschitz continuous*

(ii) $\frac{\sigma_{\min \neq 0}^2(\mathbf{C})}{\beta}$-*strongly monotone* $\forall \mathbf{z} \in \mathrm{ran}(\mathbf{C})$ *where $\sigma_{\min \neq 0}(\mathbf{C})$ is the smallest nonzero singular value of $\mathbf{C}$.*

*Proof.* Under the assumption that $f \in \Gamma_0$ is $\sigma$-strongly convex and $\beta$-smooth, from Theorem 18.15 [34], $f^*$ is both $\frac{1}{\beta}$-strongly convex and $\frac{1}{\sigma}$-smooth. It follows that $\partial f^*$ is therefore both $\frac{1}{\beta}$ strongly monotone and $\frac{1}{\sigma}$ Lipschitz continuous.

In the case of *(i)*, due to the Lipschitz continuity of $\partial f^*$

$$\begin{aligned}\|\mathbf{C}\left(\partial f^*(\mathbf{C}^T \mathbf{z}_1) - \partial f^*(\mathbf{C}^T \mathbf{z}_2)\right)\| \\ \leq \sigma_{\max}(\mathbf{C}) \|\partial f^*(\mathbf{C}^T \mathbf{z}_1) - \partial f^*(\mathbf{C}^T \mathbf{z}_2)\| \\ \leq \frac{\sigma_{\max}(\mathbf{C})}{\mu}\|\mathbf{C}^T(\mathbf{z}_1 - \mathbf{z}_2)\| \\ \leq \frac{\sigma_{\max}^2(\mathbf{C})}{\mu}\|\mathbf{z}_1 - \mathbf{z}_2\|,\end{aligned}$$

As such, $\mathbf{C}\partial f^*(\mathbf{C}^T\bullet)$ is $\frac{\sigma_{\max}(\mathbf{C})^2}{\mu}$-Lipschitz continuous.

In the case of *(ii)*, due to the strong monotonicity of $\partial \mathbf{f}^*$

$$\begin{aligned}\langle \mathbf{C}\left(\partial f^*(\mathbf{C}^T\mathbf{z}_1) - \partial f^*(\mathbf{C}^T\mathbf{z}_2)\right), \mathbf{z}_1 - \mathbf{z}_2\rangle \\ = \langle \partial f^*(\mathbf{C}^T\mathbf{z}_1) - \partial f^*(\mathbf{C}^T\mathbf{z}_1), \mathbf{C}^T(\mathbf{z}_1 - \mathbf{z}_2)\rangle \\ \geq \frac{1}{\beta}\|\mathbf{C}^T(\mathbf{z}_1 - \mathbf{z}_2)\|^2.\end{aligned}$$

$\forall \mathbf{z}_1, \mathbf{z}_2 \in \mathrm{ran}(\mathbf{C})$ it follows that

$$\frac{1}{\beta}\|\mathbf{C}^T(\mathbf{z}_1 - \mathbf{z}_2)\|^2 \geq \frac{\sigma^2_{\min\neq 0}(\mathbf{C})}{\beta}\|\mathbf{z}_1 - \mathbf{z}_2\|^2,$$

completing the proof. $\square$

**Remark 5.** *This proof reflects that of the method adopted in [45] in the centralised case. Note however that the similar result demonstrated therein does not hold in the distributed context due to the row-rank deficiency of the matrix $\mathbf{C}$. Specifically, Assumption 2 of [45] is violated.*

From [45], as $\mathbf{C}\partial f^*(\mathbf{C}^T\bullet)$ is both strongly monotone and Lipschitz continuous over $\mathrm{ran}(\mathbf{C})$, $\mathbf{R}_{\mathbf{T}_1,\rho}$ is contractive $\forall \mathbf{z} \in \mathrm{ran}(\mathbf{C})$ with an upper bound on this contraction given by

$$\delta = \max\left(\frac{\rho\frac{\sigma^2_{\max}(\mathbf{C})}{\mu} - 1}{\rho\frac{\sigma^2_{\max}(\mathbf{C})}{\mu} + 1}, \frac{1 - \rho\frac{\sigma^2_{\min\neq 0}(\mathbf{C})}{\beta}}{1 + \rho\frac{\sigma^2_{\min\neq 0}(\mathbf{C})}{\beta}}\right) \in [0,1).$$

The "partially contractive" nature of $\mathbf{R}_{\mathbf{T}_1,\rho}$ therefore leads to the following results. Let $\tilde{\mathbf{z}}_0$ be defined as per (19) and $\tilde{\mathbf{y}}_1 = \mathbf{R}_{\mathbf{T}_1,\rho}(\tilde{\mathbf{z}}_0)$. As $\mathbf{y}^{(k+1)} = \mathbf{R}_{\mathbf{T}_1,\rho}(\mathbf{z}^{(k)})$, it follows that,

$$\begin{aligned}\|\mathbf{y}^{(k+1)} - \tilde{\mathbf{y}}_1\|^2 &\leq \delta^2 \|\underset{\mathrm{ran}(\mathbf{C})}{\Pi}\left(\mathbf{z}^{(k)} - \tilde{\mathbf{z}}_0\right)\|^2 \\ &+ \|\underset{\mathrm{ker}(\mathbf{C})}{\Pi}\left(\mathbf{z}^{(k)} - \tilde{\mathbf{z}}_0\right)\|^2,\end{aligned} \quad (27)$$

where the vector

$$\underset{\mathrm{ker}(\mathbf{C})}{\Pi}\left(\mathbf{y}^{(k+1)} - \tilde{\mathbf{y}}_1\right) = \underset{\mathrm{ker}(\mathbf{C})}{\Pi}\left(\mathbf{z}^{(k)} - \tilde{\mathbf{z}}_0\right).$$

By the same arguments as used in Lemma VI.1, the operator $\mathbf{P} \circ \mathbf{R}_{\mathbf{T}_1,\rho} \circ \mathbf{P}$ is $\delta$ contractive over $\mathrm{ran}(\mathbf{PC})$ so that,

$$\begin{aligned}\|\mathbf{z}^{(k+2)} - \tilde{\mathbf{z}}_0\|^2 &\leq \delta^2 \|\underset{\mathrm{ran}(\mathbf{PC})}{\Pi}\left(\mathbf{y}^{(k+1)} - \tilde{\mathbf{y}}_1\right)\|^2 \\ &+ \|\underset{\mathrm{ker}(\mathbf{PC})}{\Pi}\left(\mathbf{y}^{(k+1)} - \tilde{\mathbf{y}}_1\right)\|^2,\end{aligned} \quad (28)$$

where

$$\underset{\mathrm{ker}(\mathbf{PC})}{\Pi}\left(\mathbf{z}^{(k+2)} - \tilde{\mathbf{z}}_0\right) = \underset{\mathrm{ker}(\mathbf{PC})}{\Pi}\left(\mathbf{y}^{(k+1)} - \tilde{\mathbf{y}}_1\right).$$

The two-step PDMM iteration given by

$$\mathbf{z}^{(k+2)} = (\mathbf{P} \circ \mathbf{R}_{\mathbf{T}_1,\rho} \circ \mathbf{P}) \circ \mathbf{R}_{\mathbf{T}_1,\rho}\left(\mathbf{z}^{(k)}\right),$$

can therefore be thought of as the composition of the operators $\mathbf{R}_{\mathbf{T}_1,\rho}$ and $\mathbf{P} \circ \mathbf{R}_{\mathbf{T}_1,\rho} \circ \mathbf{P}$ with each being $\delta$-contractive over $\mathrm{ran}(\mathbf{C})$ and $\mathrm{ran}(\mathbf{PC})$ respectively.

While the "partially contractive" nature of the PDMM updates suggests its geometric convergence, it is unclear what this convergence rate may be. For this reason, in the following we derive a geometrically primal convergence bound by connecting two-step PDMM with the GAP algorithm [44].

*C. A Geometric Rate Bound for PDMM Based on the Generalised Alternating Method of Projections*

Consider the computation of an upper bound on the geometric contraction factor of the sequence $\left(\mathbf{z}^{(2k)} - \tilde{\mathbf{z}}_0\right)_{k\in\mathbb{Z}_+}$ where $\tilde{\mathbf{z}}_0$ is defined as per (19). This bound can be computed as the solution of the optimisation problem

$$\begin{aligned}\max_{\mathbf{z}^{(k+2)}} \quad & \|\mathbf{z}^{(k+2)} - \tilde{\mathbf{z}}_0\|^2 \\ \text{s.t.} \quad & \mathbf{y}^{(k+1)} = \mathbf{R}_{\mathbf{T}_1,\rho}\left(\mathbf{z}^{(k)}\right) \\ & \mathbf{z}^{(k+2)} = \mathbf{P} \circ \mathbf{R}_{\mathbf{T}_1,\rho} \circ \mathbf{P}\left(\mathbf{y}^{(k+1)}\right) \\ & \|\mathbf{z}^{(k)} - \tilde{\mathbf{z}}_0\|^2 \leq 1.\end{aligned} \quad (29)$$

For this particular choice of $\tilde{\mathbf{z}}_0$, as in (20), it follows that

$$\|\mathbf{z}^{(k+2)} - \tilde{\mathbf{z}}_0\|^2 = \|\underset{\mathrm{ran}(\mathbf{PC})+\mathrm{ran}(\mathbf{C})}{\Pi}\left(\mathbf{z}^{(k+2)} - \tilde{\mathbf{z}}_0\right)\|^2.$$

Therefore, applying, (27) and (28), (29) is equivalent to

$$\max_{\mathbf{z}^{(k)}} \quad \|\left(\delta\underset{\mathrm{ran}(\mathbf{PC})}{\Pi} + \underset{\mathrm{ker}(\mathbf{PC})}{\Pi}\right)\left(\mathbf{y}^{(k+1)} - \tilde{\mathbf{y}}_1\right)\|^2$$

$$\text{s.t.} \quad \|\underset{\mathrm{ran}(\mathbf{C})}{\Pi}\left(\mathbf{y}^{(k+1)} - \tilde{\mathbf{y}}_1\right)\|^2 \leq \delta^2 \|\underset{\mathrm{ran}(\mathbf{C})}{\Pi}\left(\mathbf{z}^{(k)} - \tilde{\mathbf{z}}_0\right)\|^2 \quad (30\mathrm{a})$$

$$\underset{\mathrm{ker}(\mathbf{C})}{\Pi}\left(\mathbf{y}^{(k+1)} - \tilde{\mathbf{y}}_1\right) = \underset{\mathrm{ker}(\mathbf{C})}{\Pi}\left(\mathbf{z}^{(k)} - \tilde{\mathbf{z}}_0\right) \quad (30\mathrm{b})$$

$$\|\mathbf{z}^{(k)} - \tilde{\mathbf{z}}_0\|^2 \leq 1, \quad (30\mathrm{c})$$

where in the objective we have exploited the orthogonality of $\mathrm{ran}(\mathbf{PC})$ and $\mathrm{ker}(\mathbf{PC})$. The constraints (30a), (30b) and (30c) collectively define the feasible set of the vectors $\mathbf{y}^{(k+1)} - \tilde{\mathbf{y}}_1$. As the set of the feasible vectors $\mathbf{z}^{(k)} - \tilde{\mathbf{z}}_0$ (30c) denotes a sphere, the set of $\mathbf{y}^{(k+1)} - \tilde{\mathbf{y}}_1$ is an ellipsoid given by

$$\mathcal{E}_\mathbf{y} = \left\{\left(\delta\underset{\mathrm{ran}(\mathbf{C})}{\Pi} + \underset{\mathrm{ker}(\mathbf{C})}{\Pi}\right)\mathbf{u} \mid \|\mathbf{u}\| \leq 1\right\}.$$

The optimisation problem in (29) is therefore equivalent to

$$\begin{aligned}\max_{\mathbf{z}^{(k)}} \quad & \|\left(\underset{\mathrm{ker}(\mathbf{PC})}{\Pi} + \delta\underset{\mathrm{ran}(\mathbf{PC})}{\Pi}\right)\left(\underset{\mathrm{ker}(\mathbf{C})}{\Pi} + \delta\underset{\mathrm{ran}(\mathbf{C})}{\Pi}\right)\mathbf{u}\|^2 \\ \text{s.t.} \quad & \|\mathbf{u}\|^2 \leq 1 \\ & \mathbf{u} \in \mathrm{ran}(\mathbf{PC}) + \mathrm{ran}(\mathbf{C}),\end{aligned}$$

where the additional domain constraint stems from (20). Furthermore, this is equivalent to

$$\begin{aligned}\max_{\mathbf{z}^{(k)}} \quad & \|\left((1-\delta)\underset{\mathrm{ker}(\mathbf{PC})}{\Pi} + \delta\mathbf{I}\right)\left((1-\delta)\underset{\mathrm{ker}(\mathbf{C})}{\Pi} + \delta\mathbf{I}\right)\mathbf{u}\|^2 \\ \text{s.t.} \quad & \|\mathbf{u}\|^2 \leq 1 \\ & \mathbf{u} \in \mathrm{ran}(\mathbf{PC}) + \mathrm{ran}(\mathbf{C}),\end{aligned} \quad (31)$$

which corresponds to finding the convergence rate bound for a particular instance of the GAP algorithm [44]. In the following, we will denote the GAP operator via

$$\hat{\mathbf{T}}^2_{P,\rho} = \left((1-\delta)\underset{\mathrm{ker}(\mathbf{PC})}{\Pi} + \delta\mathbf{I}\right)\left((1-\delta)\underset{\mathrm{ker}(\mathbf{C})}{\Pi} + \delta\mathbf{I}\right).$$

Note that, as $\hat{\mathbf{T}}^2_{P,\rho}$ is the composition of nonexpansive operators and is similar to the positive semidefinite matrix,

$$\mathcal{T} = \left(\prod_{\ker(\mathbf{PC})} + \delta \prod_{\operatorname{ran}(\mathbf{PC})}\right)^{-\frac{1}{2}} \hat{\mathbf{T}}^2_{P,\rho} \left(\prod_{\ker(\mathbf{PC})} + \delta \prod_{\operatorname{ran}(\mathbf{PC})}\right)^{\frac{1}{2}},$$

the spectrum $\sigma\left(\hat{\mathbf{T}}^2_{P,\rho}\right) \in [0,1]$. Furthermore, $\forall \lambda \in \sigma\left(\hat{\mathbf{T}}^2_{P,\rho}\right) \mid \lambda = 1$, by inspection the corresponding eigenvector $\mathbf{v} \in \ker(\mathbf{C}) \cap \ker(\mathbf{PC})$. In the case of (31), the convergence rate will therefore be determined by the subdominant eigenvalue of $\hat{\mathbf{T}}^2_{P,\rho}$, given by

$$\gamma = \max\left\{\lambda \mid \lambda \in \sigma\left(\hat{\mathbf{T}}^2_{P,\rho}\right), \lambda \neq 1\right\} < 1.$$

Optimal rate bounds for generalisations of the classic alternating projections algorithm has been an area of recent attention in the literature with two notable papers on the subject being [46] and [44]. Specifically, determining $\gamma$ for the GAP problem was considered in [44]. In particular, applying basic algebraic manipulation to equation 8 of [44] leads to

$$\gamma = \left|\max\left\{\delta + (1-\delta)\cos(\theta_F)\left(\frac{(1-\delta)\cos(\theta_F)}{2}\right.\right.\right. \\ \left.\left.\left.\pm \sqrt{\frac{(1-\delta)^2\cos^2(\theta_F)}{4} + \delta}\right)\right\}\right| < 1, \quad (32)$$

where $\theta_F > 0$ denotes the Freidrichs angle, the smallest nonzero principle angle, between the subspaces $\operatorname{ran}(\mathbf{C})$ and $\operatorname{ran}(\mathbf{PC})$. It follows that, for two-step PDMM,

$$\|\mathbf{z}^{(k+2)} - \tilde{\mathbf{z}}_0\|^2 \leq \gamma^2 \|\mathbf{z}^{(k)} - \tilde{\mathbf{z}}_0\|^2.$$

For even $k$, we therefore obtain

$$\|\mathbf{z}^{(k+1)} - \mathbf{T}_{P,\rho}(\tilde{\mathbf{z}}_0)\|^2 \leq \gamma^k \|\mathbf{z}^1 - \tilde{\mathbf{z}}_1\|^2 \\ \leq \gamma^k \|\mathbf{z}^0 - \tilde{\mathbf{z}}_0\|^2,$$

so that a geometrically converging upper bound for the auxiliary error of PDMM is given by

$$\|\mathbf{z}^{(k+2)} - \tilde{\mathbf{z}}_0\|^2 \leq \gamma^{k+2} \frac{\|\mathbf{z}^0 - \tilde{\mathbf{z}}_0\|^2}{\gamma}.$$

Applying the inequality in (24) it follows that

$$\|\mathbf{x}^{(k+2)} - \mathbf{x}^*\|^2 \leq \gamma^{k+2} \frac{\sigma^2_{\max}(\mathbf{C})}{\rho \sigma^2_{\min\neq 0}(\mathbf{C})} \frac{\|\mathbf{z}^{(0)} + \tilde{\mathbf{z}}_0\|^2}{\gamma} \\ = \gamma^{k+2} \epsilon$$

Thus $\mathbf{x}^{(k+2)} \to \mathbf{x}^*$ is upper bounded by a geometrically contracting sequence which, to the best of the authors knowledge, is a first for PDMM within the literature.

### D. Step Size Selection for PDMM

We now move to demonstrating how we can optimise the geometric convergence rate bound introduced in Sec. VI to improve the performance guarantees of PDMM. In particular, we consider the optimal selection of the step size $\rho$ for a given problem instance. The proposed method is inherently distributed and only requires a finite number of communications between nodes, making it practical to implement.

For a given instance of problem (1), with $\theta_F$ fixed (but not nescessarily known), $\gamma$, as defined in (32) is monotonically increasing with the contraction factor $\delta$. Thus, convergence rate can be maximised by selecting an appropriate step size $\rho$ to minimize $\delta$. In [45] this optimal $\rho$ was shown to be

$$\rho^* = \frac{\sqrt{\beta\mu}}{\sigma_{\max}(\mathbf{C})\sigma_{\min\neq 0}(\mathbf{C})}, \quad (33)$$

resulting in a corresponding subspace contraction factor

$$\delta^* = \frac{\sqrt{\kappa}-1}{\sqrt{\kappa}+1}, \ \kappa = \frac{\sigma^2_{\max}(\mathbf{C})\beta}{\sigma^2_{\min\neq 0}(\mathbf{C})\mu}.$$

The problem of optimal step size selection therefore reduces to determining the largest and smallest non-zero singular values of the matrix $\mathbf{C}$ as well as the strong convexity ($\mu$) and smoothness ($\beta$) constants of the function $f$. The latter two parameters can be computed as the minimum $\mu_i$ and maximum $\beta_i$ across all nodes. Furthermore, from the block diagonal structure of the matrix $\mathbf{C}$, the singular values can also be computed in a parallel fashion across nodes. In particular, each node $i$ can determine its local maximum and minimum (non-zero) singular values which we will denote by $\sigma_{\max,i} = \max\{\sigma(\mathbf{C}_i)\}$ and $\sigma_{\min\neq 0,i} = \min\{\sigma(\mathbf{C}_i) > 0\}$ respectively. The task, therefore, is to find the variables

$$\sigma_{\max}(\mathbf{C}) = \max\{\sigma_{\max,i}, \forall i \in V\} \\ \sigma_{\min\neq 0}(\mathbf{C}) = \min\{\sigma_{\min\neq 0,i}, \forall i \in V\} \\ \mu = \max\{\mu_i, \forall i \in V\} \\ \beta = \min\{\beta_i, \forall i \in V\},$$

so that $\rho^*$ can be computed as per (33).

In the context of a distributed network, intuitively each node can learn the maximum or minimum of a set by repeatedly exchanging local estimates of these values with their neighbours. This allows the extreme values of a set to be diffused into a network regardless of topology. Importantly, this process requires only $D$ transmissions per node, where $D$ is the diameter of the underlying network, allowing each node to compute the optimal step size choice within a finite number of transmissions. This method is summarised in Algorithm 4. Note that this process is only executed once, at the instantiation of the PDMM problem and therefore comes at a negligible transmission cost in contrast to the iteration of the algorithm.

## VII. NUMERICAL EXPERIMENTS

In this section, we verify the analytical results of Sec. V and VI with numerical experiments. These results are broken down into three subsections.

### A. PDMM for Differentiable and/or Strictly Convex Functions

The first set of simulations aims to validate the sufficiency of strict convexity and differentiability to guarantee primal convergence as introduced in Sec. V-B. For these simulations, as testing all such functions would be computationally infeasible, we instead considered the family $p$-th power of $p$-norms

**Algorithm 4** Optimal Step Size Selection for PDMM
1: **Initialise:**
2: $\sigma_{\max,i}^{(0)}(\mathbf{C}_j) = \max\{\sigma(\mathbf{C}_i)\} \ \forall i \in V$
3: $\sigma_{\min \neq 0,i}^{(0)}(\mathbf{C}_i) = \min\{\sigma(\mathbf{C}_i) > 0\} \ \forall i \in V$
4: $\mu_i^{(0)}, \beta_i^{(0)} \in \mathbb{R} \ \forall i \in V$
5: **for** k=0,...,D-1 **do**    ▷ Iterative Set Max and Min
6:     **for all** $i \in V$ **do**
7:         $\sigma_{\max,i}^{(k+1)} = \max\{\sigma_{\max,j}^{(k)} \forall j \in \mathcal{N}(i) \cup i\}$
8:         $\sigma_{\min \neq 0,i}^{(k+1)} = \max\{\sigma_{\min \neq 0,j}^{(k)} \forall j \in \mathcal{N}(i) \cup i\}$
9:         $\mu_i^{(k+1)} = \min\{\mu_j^{(k)} \forall j \in \mathcal{N}(i) \cup i\}$
10:       $\beta_i^{(k+1)} = \max\{\beta_j^{(k)} \forall j \in \mathcal{N}(i) \cup i\}$
11:     **end for**
12: **end for**
13: **for all** $i \in V$ **do**    ▷ Compute Optimal Step Size
14:     $\rho_i^* = \dfrac{\sqrt{\beta_i^{(D)} \mu_i^{(D)}}}{\sigma_{\max,i}^{(D)} \sigma_{\min \neq 0,i}^{(D)}}$
15: **end for**

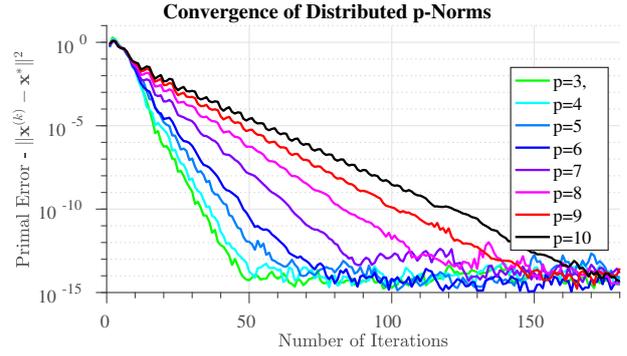

Fig. 1: The primal convergence of different $p$-norm$^p$ consensus problem for a 10 node Erdős-Rényi network.

for $p \in \{3, 4, 5, \cdots\}$. The prototype optimisation problem for these simulations was given by

$$\min_{\mathbf{x}} \sum_{i \in V} \|\mathbf{x}_i - \mathbf{a}_i\|_p^p$$
$$\text{s.t.} \quad \mathbf{x}_i - \mathbf{x}_j = \mathbf{0} \quad \forall (i,j) \in E,$$

where $\mathbf{a}_i$ are local observation vectors and, for simplicity, edge based consensus constraints were chosen.

An $N = 10$ node undirected Erdős-Rényi network [47] was considered for these simulations. Such a network is based on a random graph model where $\forall \ i, \ j \in V, i \neq j$, there is equal probability that $(i,j) \in E$. This probability determines the density of the connectivity in the network and in this case was set to $\frac{\log(N)}{N}$. This resulted in a total of 12 undirected edges and a resulting network which was verified as forming a single connected component as per the assumptions in Sec. IV. The resulting $\mathbf{C}$ and $\mathbf{P}$ matrices were generated as per (4) and were used for all simulations. Additionally, a randomly generated initial auxiliary vector $\mathbf{z}^{(0)}$ was also used for all problem instances.

For $p = 3, \cdots 10$, 180 iterations of PDMM were performed and the resulting primal error computed. Figure 1 demonstrates the convergence of this error with respect to iteration count. The squared Euclidean distance between the primal iterates and the primal optimal set was used as a measure of convergence for these simulations. Such a measure has the additional benefit of being stricter than all of the $p$-norms considered and thus also guarantees convergence in these metrics. As the purpose of this figure was to verify the convergence guarantee rather than demonstrate convergence rate, the step size $\rho$ was selected for each instance to improve visual clarity. Note that the final finite precision stems from the use of MATLABs $fminunc$ function during each primal update.

Figure 2 further demonstrates that the choice of $\rho$ does not effect the guarantee of convergence. In this instance, the number of iterations required to reach an auxiliary precision of $1e^{-5}$ was used as a surrogate for overall algorithmic convergence. This measure was chosen as the auxiliary error of PDMM forms a monotonically decreasing sequence. In contrast the primal error need not satisfy this, a point which can be verified in Figure 1. Note that while there is a clear variation in the rate of convergence for different choices of $\rho$, the guarantee of convergence of the algorithms is unaffected.

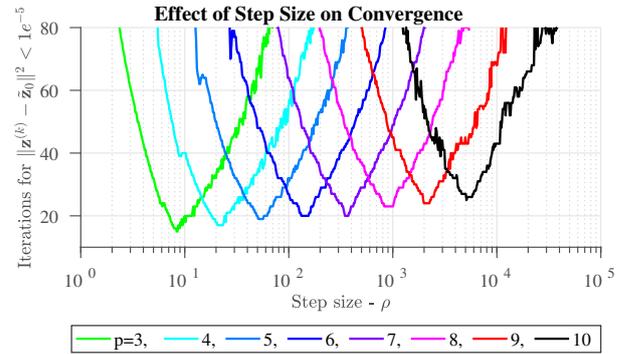

Fig. 2: A comparison to the required iterations for $\|\mathbf{z}^{(k)} - \tilde{\mathbf{z}}_0\|^2 \leq 1e^{-5}$ for variation in step size ($\rho$).

### B. Averaged PDMM for General Functions in $\Gamma^0$

We now move our attention to demonstrating, via contradiction, that $f \in \Gamma_0$ alone is insufficient to guarantee the convergence of PDMM without additional modification such as operator averaging. We highlight this point as while it is well known with respect to PR splitting, currently within the literature, there are no results suggesting this point for PDMM.

For the purpose of this example we will consider the problem of distributed L1 consensus, given by

$$\min_{\mathbf{x}} \sum_{i \in V} \|\mathbf{x}_i - \mathbf{a}_i\|_1$$
$$\text{s.t.} \quad \mathbf{x}_i - \mathbf{x}_j = \mathbf{0} \quad \forall (i,j) \in E.$$

To demonstrate the contradiction, we need only demonstrate a single problem instance for which convergence is not achieved.

For this purpose the same network from section VII-A was used and a step size of $\rho = \frac{1}{2}$ was selected. Figure 3 demonstrates the observed convergence plots of both standard PDMM and ADMM for this problem. Due to the lack of uniqueness of the optimal variables, convergence was measured via objective sub-optimality given by $\|f(\mathbf{x}^{(k)}) - f(\mathbf{x}^*)\|^2$.

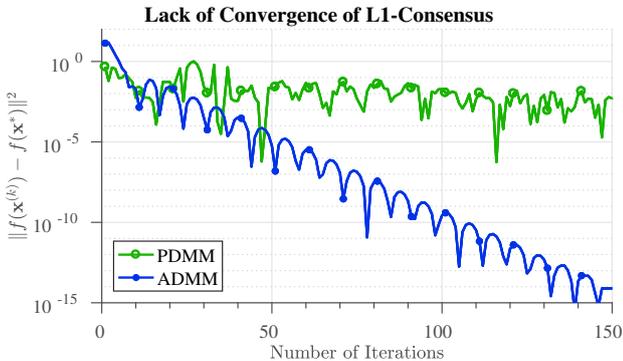

Fig. 3: A comparison of the (lack of) convergence of PDMM (green) and averaged PDMM (blue) based L1-Consensus. Without the use of averaging primal convergence cannot be guaranteed for PDMM for general functions $f \in \Gamma_0$.

One can clearly observe the lack of convergence of PDMM for the considered L1 problem with the algorithm instead oscillating around a suboptimal value. In contrast, the averaged PDMM instance does not suffer from this shortcoming. As such, $f \in \Gamma_0$ is not sufficient to guarantee primal optimality for general PDMM functions.

### C. Geometric Convergence of PDMM for Strongly Convex and Smooth Functions

Finally we move to verifying the geometric convergence bound proposed in Sec. VI. We do so by comparing the convergence of multiple problem instances to this bound. For these simulations, 10000 Erdős-Rényi networks, each comprised of 10 nodes, were generated with a connection probability of $\frac{\log(N)}{N}$. The considered network instances were all verified as forming single connected networks. Each network was then paired with a randomly generated distributed quadratic objective function. The local variables at each node were chosen to be three element vectors so that $\mathbf{x}_i \in \mathbb{R}^3 \ \forall i \in V$ resulting in a total of 30 variables. The constraint functions were again chosen to be consensus constraints between neighbouring nodes. The resulting problem instances were therefore all of the form

$$\min_{\mathbf{x}} \sum_{i \in V} \left(\frac{1}{2}\mathbf{x}_i^T \mathbf{Q}_i \mathbf{x}_i - \mathbf{q}_i^T \mathbf{x}_i\right)$$
$$\text{s.t.} \quad \mathbf{x}_i - \mathbf{x}_j = \mathbf{0} \quad \forall (i,j) \in E.$$

For each problem instance, the matrices $\mathbf{Q}_i \succeq \mathbf{0}$ were generated in such a way that a constant convergence rate bound was achieved. In this case the contraction factor of this rate bound was specified as $\gamma = 0.9$. Furthermore, the initial vectors $\mathbf{z}^{(0)}$ were generated randomly and for each the associated $\tilde{\mathbf{z}}_0$ was computed as per Eq. (19). This randomisation procedure was implemented so that $\epsilon = 1$ for all instances. Thus, a single global convergence bound was considered for all realisations.

For each problem instance a total of 120 iterations of PDMM, as per Algorithm 3, were performed and the auxiliary errors, $\|\mathbf{z}^{(k)} - \tilde{\mathbf{z}}_0\|^2$ for $k$ even and $\|\mathbf{z}^{(k)} - \tilde{\mathbf{z}}_1\|^2$ for $k$ odd, were computed for each iteration. The distribution of the resulting data is demonstrated in Figure 4 which highlights both the worst case convergence and the spread of the convergence curves for all problem instances.

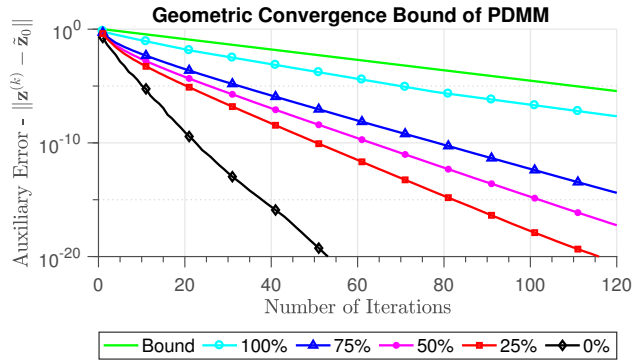

Fig. 4: Convergence of simulated PDMM problem instances. From top to bottom, the solid green line denotes the convergence rate bound while the remaining 5 lines denote the 100%, 75%, 50%, 25% and 0% quantiles respectively.

As expected, the proposed bound holds for all instances. One can further note that the rate of the worst case sequence (100% quantile) does not exceed that of the bound. Additionally, while the bound holds for the worst case functions, most problem instances exhibit far faster convergence than that predicted. This suggests that for more restrictive problem classes stronger bounds may be able to be guaranteed. While the $\rho$ selection defined in Sec. VI-D may optimise the proposed convergence rate bound, it is highly unlikely that it is optimal for all problems and thus a more refined parameter selection process may be needed for a given problem instance.

## VIII. CONCLUSIONS

In this paper we have presented a novel derivation of the node-based distributed algorithm termed the primal-dual method of multipliers (PDMM). Unlike existing efforts within the literature, monotone operator theory was used for this purpose, providing both a succinct derivation for the algorithm while highlighting the relationship between PDMM and other existing first order methods such as PR splitting and ADMM. Furthermore, the proposed derivation lead to a collection of new results for the algorithm. Guaranteed primal convergence was demonstrated for strictly convex, differentiable functions and, in the case of strongly convex, smooth local functions, a geometric convergence bound was presented, a first for PDMM. In conclusion the demonstrated results unify PDMM with existing solvers in the literature while providing new insight into its operation and convergence characteristics.